\documentclass[2p]{article}
\usepackage{amssymb,amsmath,amsthm}
\usepackage{indentfirst}
\usepackage{exscale}
\usepackage{relsize}
\usepackage[numbers,sort&compress]{natbib}

\usepackage{geometry}
\usepackage{color}
\newcommand{\R}{\mathbb{R}}
\theoremstyle{definition}

\allowdisplaybreaks
\theoremstyle{remark}

\numberwithin{equation}{section}
 \UseRawInputEncoding

\begin{document}
\title{\Large\bf{ Dependence on parameters of solutions for a generalized poly-Laplacian system on weighted graphs}}
\date{}
\author {Xiaoyu Wang$^{1}$, Junping Xie$^{2}$ \footnote{Corresponding author, E-mail address: hnxiejunping@163.com}, \ Xingyong Zhang$^{1,3}$, \ Xin Ou$^{1}$\\
{\footnotesize $^1$Faculty of Science, Kunming University of Science and Technology, Kunming, Yunnan, 650500, P.R. China.}\\
{\footnotesize $^2$Faculty of Transportation Engineering, Kunming University of Science and Technology,}\\
 {\footnotesize Kunming, Yunnan, 650500, P.R. China.}\\
{\footnotesize $^{3}$Research Center for Mathematics and Interdisciplinary Sciences, Kunming University of Science and Technology,}\\
 {\footnotesize Kunming, Yunnan, 650500, P.R. China.}\\}

 \date{}
 \maketitle

 \begin{center}
 \begin{minipage}{15cm}
 \par
 \small  {\bf Abstract:}  We mainly investigate the continuous dependence on parameters of nontrivial  solutions for a generalized poly-Laplacian system on the weighted finite graph $G=(V, E)$. We firstly present an existence result of mountain pass type nontrivial solutions when the nonlinear term $F$ satisfies the super-$(p, q)$ linear growth condition which is a simple generalization of those results in \cite{Zhang2022}. Then we mainly show that the mountain pass type nontrivial solutions of the poly-Laplacian system are uniformly bounded for parameters and  the concrete upper and lower bounds are given, and are continuously dependent on parameters. Similarly, we also present the existence result, the concrete upper and lower bounds, uniqueness, and  dependence on parameters for the locally minimum type nontrivial solutions.  Subsequently, we present an example on optimal control as an application of our results. Finally, we give a nonexistence result and some results for the  corresponding scalar equation.
 \end{minipage}
 \end{center}
  \allowdisplaybreaks
 \vskip2mm
 {\section{Introduction }}
\setcounter{equation}{0}
The main aim of this paper is to study the continuous dependence on parameters of nontrivial solutions for  the following generalized poly-Laplacian system on weighted finite graph $G=(V,E)$:
\begin{eqnarray}
\label{eq1}
 \begin{cases}
  \pounds_{m_1,p}u+h_1(x)|u|^{p-2}u=F_u(x,u,v,w),\;\;\;\;\hfill x\in V,\\
  \pounds_{m_2,q}v+h_2(x)|v|^{q-2}v=F_v(x,u,v,w),\;\;\;\;\hfill x\in V,\\
   \end{cases}
\end{eqnarray}
where  $V$ is the vertexes set of finite graph $G$, $E$ is the edges set of finite graph $G$, $m_i\in \mathbb{N}, i=1,2$, and $p,q\ge 2$, $h_i:V\to \R^+,i=1,2$, and $F:V\times \R^2\times \R \to \R$, $w\in \mathcal{J}$, $\mathcal{J} \subset\R$ is a bounded closed interval, and $\pounds_{m,s}, (m=m_1,m_2, s=p,q)$ are  defined as following: for any function  $\phi:V\to\R$,
\begin{eqnarray}\label{eq2}
\int_V (\pounds_{m,s}u)\phi d\mu=
 \begin{cases}
  \int_V|\nabla^m u|^{s-2}\Gamma(\Delta^{\frac{m-1}{2}}u,\Delta^{\frac{m-1}{2}}\phi)d\mu,& \text { when $m$ is odd},\\
  \int_V|\nabla^m u|^{s-2}\Delta^{\frac{m}{2}}u\Delta^{\frac{m}{2}}\phi d\mu,&  \text { when $m$ is even}.
    \end{cases}
\end{eqnarray}
When $s=2$, $\pounds_{m,s}=(-\Delta)^mu$ is referred to as the poly-Laplacian operator of $u$ and when $m=1$ and $s=p$, $\pounds_{m,p}=-\Delta_pu$, which is called $p$-Laplacian operator defined by (\ref{Eq6}) below.
\par
When $m_1=m_2=m,p=q$, $u=v$ and $h_1=h_2=h$, system (\ref{eq1}) degenerates into the following equation:
\begin{eqnarray}
\label{eq3}
  \pounds_{m,p}u+h(x)|u|^{p-2}u=f(x,u,w),\;\;\;\;\hfill x\in V,
\end{eqnarray}
where $F(x,u,w)=\int_0^u f(x,u,w)du$ for all $x\in V$ and $w\in \mathcal{J}$. Corresponding to system (\ref{eq1}), we also present a result for equation (\ref{eq3}) in section 7 below.
\par
In recent years, the research on the existence and multiplicity of solutions for partial differential equations on graphs has attracted some attentions, see \cite{Ge20181, Ge20184, Grigor yan2016, Han2021, Yu2024, Zhang2022, Zhang2018, Zhang2019} for related topics. In particular, in \cite{Grigor yan2016}, Grigor'yan-Lin-Yang used the mountain pass theorem to obtain the existence of a positive solution for Yamabe equation, $p$-Laplacian equation and (\ref{eq3}) under the well-known Ambrosetti-Rabinowitz condition (see \cite{Ambrosetti1973}). It is worth noting that the basic variational framework and some embedding theorems are given in \cite{Grigor yan2016}, and the Sobolev embedding inequality on graphs is developed, which lays a foundation for the application of the minimax method \cite{Chang1986,Rabinowitz1986} to the partial differential equations on weighted graphs.
\par
Inspired by \cite{Grigor yan2016}, Zhang-Zhang-Xie-Yu \cite{Zhang2022} studied the following high-order Yamabe-type coupled system on finite graphs:
\begin{eqnarray}
\label{eq4}
 \begin{cases}
  \pounds_{m_1,p}u+h_1(x)|u|^{p-2}u=F_u(x,u,v),\;\;\;\;\hfill x\in V,\\
  \pounds_{m_2,q}v+h_2(x)|v|^{q-2}v=F_v(x,u,v),\;\;\;\;\hfill x\in V.
   \end{cases}
\end{eqnarray}
When the nonlinear term satisfies the super-$(p,q)$ linear growth condition, the existence and multiplicity of the nontrivial solutions for system (\ref{eq4}) are obtained by using the mountain pass theorem and the symmetric mountain pass theorem, respectively. Yu-Zhang-Xie-Zhang \cite{Yu2024} studied the existence of nontrivial solutions for a class of generalized poly-Laplacian nonlinear system with mixed nonlinear terms on finite graphs and a generalized  poly-Laplacian nonlinear system on locally finite graphs with Dirichlet boundary by the mountain pass theorem, when $F$ satisfies the asymptotically quadratic condition at infinity. We refer readers to \cite{Yang2023,Yang2024,Pang2023,Pang2024} for more results on the generalized poly-Laplacian system and $(p,q)$-Laplacian system.
\par
There is another case of discrete-like partial differential equations other than the one on the graphs, which is the fractals case, and some works have been done, for instance, see \cite{Bonanno2012, Molica Bisci2017, Breckner2010, Breckner2011, Galewski2011, Galewski2019, Falconer1999}. Especially, Marek Galewski \cite{Galewski2019} considered the second-order elliptic equation with a Dirichlet boundary condition on a Sierpi\'{n}ski gasket which is a fractal domain:
\begin{eqnarray}
\label{eq6}
 \begin{cases}
  \Delta u(x)+a(x)u(x)=f(x,u(x),w(x)),\;\;\;\;\hfill a.e.\;\; x\in V\setminus V_0,\\
  u|_{V_0}=0,\\
   \end{cases}
\end{eqnarray}
where $V$ denotes the Sierpi\'{n}ski gasket in $\R^{N-1},N\ge 2$, and $V_0$ is called the intrinsic boundary. When $f$ satisfies the superquadratic growth condition, they obtained the equation (\ref{eq6}) has a nontrivial solution by using the mountain pass theorem and  the iterative technique. They also studied the continuous parameters dependence on the mountain pass solution of the equation (\ref{eq6}). In \cite{Galewski2013} and \cite{Galewski2012}, by using the similar methods, the parameters dependence of the mountain pass solution for the equation or system with the discrete boundary value problem is also considered.
\par
Motivated by \cite{Zhang2022, Galewski2019}, in this paper, we mainly focus on the parameters dependence on nontrivial mountain pass solutions of system (\ref{eq1}). We will follow the line in \cite{Galewski2019} to discuss the parameters dependence. Moreover, as a byproduct, we obtain that the concrete range of the nontrivial solutions for each parameter $w\in \mathcal{J}$, and the range is uniform for all $w$. Smilarly, we also present the existence result, the concrete upper and lower bounds, uniqueness, and  dependence on parameters for the locally minimum type nontrivial solutions, which were not given in \cite{Galewski2019,Galewski2013,Galewski2012}.  Subsequently, we present an example on optimal control as an application of our results. Finally, we give a nonexistence result and some results for the  corresponding scalar equation.

\vskip2mm
{\section{ Preliminaries}}
\setcounter{equation}{0}
\par
In this section, we recall some notations and results. One can see more details in \cite{Grigor yan2016, Zhang2022}. Let $G=(V, E)$ be a finite graph which means both $V$ and $E$ are finite sets. For any edge $xy\in E$ connecting $x\in V$ with $y\in V$, suppose its weight satisfies $\omega_{xy}=\omega_{yx}>0$.  For any $x\in V$, define its degree as deg$(x)=\sum\limits_{y\thicksim x}\omega_{xy}$, where $y\thicksim x$ represents those $y$ connected with $x$. Suppose that $\mu:V\rightarrow \R^+$ is a finite measure. Now, we recall the definition of Laplacian operator on the graph,
\begin{eqnarray}
\label{Eq1}
\Delta u(x)=\frac{1}{\mu(x)}\sum\limits_{y\thicksim x}w_{xy}(u(y)-u(x)),\;\; \forall \; u:V\to\mathbb{R},
\end{eqnarray}
and its corresponding gradient form reads as
\begin{eqnarray}
\label{Eq2}
\Gamma(u,v)(x)=\frac{1}{2\mu(x)}\sum\limits_{y\thicksim x}w_{xy}(u(y)-u(x))(v(y)-v(x)).
\end{eqnarray}
The modulus of the gradient is defined as
\begin{eqnarray}
\label{Eq3}
|\nabla u|(x)=\sqrt{\Gamma(u,u)(x)}=\left(\frac{1}{2\mu(x)}\sum\limits_{y\thicksim x}w_{xy}(u(y)-u(x))^2\right)^{\frac{1}{2}},
\end{eqnarray}
 and the modulus of $m$-order gradient ($m\in \mathbb N$) of $u$ is defined as
\begin{eqnarray}
\label{Eq4}
|\nabla^m u|=
 \begin{cases}
  |\nabla\Delta^{\frac{m-1}{2}}u|,& \text {when $m$ is odd,}\\
  |\Delta^{\frac{m}{2}}u|,&  \text {when $m$ is even.}
   \end{cases}
\end{eqnarray}
For any function $u:V\rightarrow \R$, we denote
\begin{eqnarray*}
\label{Eq5}
\int_V u(x) d\mu=\sum\limits_{x\in V}\mu(x)u(x),
\end{eqnarray*}
and $|V|=\sum\limits_{x\in V}\mu(x)$.
\par
When $p\geq2$, the $p$-Laplacian operator $\Delta_p u$ is defined by
\begin{eqnarray}
\label{Eq6}
\Delta_p u(x)=\frac{1}{2\mu(x)}\sum\limits_{y\sim x}\left(|\nabla u|^{p-2}(y)+|\nabla u|^{p-2}(x)\right)\omega_{xy}(u(y)-u(x)).
\end{eqnarray}
In the distributive sense, we can write $\Delta_p u$ in the following form:
\begin{eqnarray*}
\label{Eq7}
\int_V(\Delta_p u)\phi d\mu=-\int_V|\nabla u|^{p-2}\Gamma(u,\phi)d\mu, \;\; \phi \in\mathcal{C}_c(V),
\end{eqnarray*}
where $\mathcal{C}_c(V)$ represents the set of all real valued functions on $V$ with compact support.
\par
Define the space
\begin{eqnarray*}
W^{m_i,s}(V)=\{u:V\to \R\},\  \ s=p,  q, i=1,2,
\end{eqnarray*}
with  the norm
\begin{eqnarray}
\label{Eq8}
\|u\|_{W^{m_i,s}(V)}=\left(\int_V(|\nabla^m u(x)|^s+h_i(x)|u(x)|^s)d\mu\right)^\frac{1}{s},  s=p,  q, i=1,2.
\end{eqnarray}
 Then $W^{m_i,s}(V)$ is a finitely dimensional Banach space. For any fixed $1<r<+\infty$, let
$$
L^r(V)=\{u:V\to \R\}
$$
endowed with the norm
\begin{eqnarray*}
\label{Eq9}
\|u\|_{L^r(V)}=\left(\int_V|u(x)|^rd\mu\right)^\frac{1}{r}.
\end{eqnarray*}

\vskip2mm
\noindent
{\bf Lemma 2.1.} (\cite{Zhang2022}) {\it For all $u\in W^{m_1,p}(V)$ and $v\in W^{m_2,q}(V)$, there hold
$$\|u\|_{\infty}\leq b\|u\|_{W^{m_1,p}(V)}, \|v\|_{\infty}\leq d\|v\|_{W^{m_2,q}(V)},$$
where $\|u\|_{\infty}=\max_{x\in V}|u(x)|, \|v\|_{\infty}=\max_{x\in V}|v(x)|$, $b=\left(\frac{1}{\mu_{\min}h_{1,\min}}\right)^\frac{1}{p}$, $d=\left(\frac{1}{\mu_{\min}h_{2,\min}}\right)^\frac{1}{q}$, $ h_{i,\min}:=\min_{x\in V}h_i(x), i=1,2, \mu_{\min}:=\min_{x\in V}\mu(x)$. }

\vskip2mm
\noindent
{\bf Lemma 2.2.} (Ekeland's variational principle \cite{Mawhin1989}) {\it Let $M$ be a  complete metric space with metric $d$, and $\varphi:M\rightarrow \R$ be a lower semicontinuous function, bounded from below and not identical to $+\infty$. Let $\varepsilon >0$ be given and $x \in M$ such that
$$
\varphi(x)\le \inf_{M} \varphi+\varepsilon.
$$
Then there exists $y \in M$ such that
$$
\varphi(y)\le \varphi(x),\;d(x,y)\le1,
$$
and for each $z\in M$, one has
$$
\varphi(y)\le \varphi(z)+\varepsilon d(y,z).
$$}

\vskip2mm
{\section{Mountain pass type solutions}}
\vskip2mm
\par
Note that the space $W:=W^{m_1,p}(V)\times W^{m_2,q}(V)$ with the norm $\|(u,v)\|=\|u\|_{W^{m_1,p}(V)}+\|v\|_{W^{m_2,q}(V)}$ is a Banach space of finite dimension. For any fixed $w\in W$, define the functional $\varphi_w:W\to \R$ as
\begin{eqnarray*}
\label{Eq11}
\varphi_w(u,v)=\frac{1}{p}\int_V(|\nabla^{m_1}u|^p+h_1(x)|u|^p)d\mu+\frac{1}{q}
\int_V(|\nabla^{m_2}v|^q+h_2(x)|v|^q)d\mu-\int_V F(x,u,v,w)d\mu.
\end{eqnarray*}
Then $\varphi_w\in C^1(W, \R)$ and for any $(u,v),(\phi_1,\phi_2)\in W$, there holds
\begin{eqnarray*}
\label{Eq12}
\langle\varphi_w'(u,v),(\phi_1,\phi_2)\rangle=\langle\varphi_{w,u}(u,v),\phi_1\rangle+\langle\varphi_{w,v}(u,v),\phi_2\rangle,
\end{eqnarray*}
where
$$
\langle\varphi_{w,u}(u,v),\phi_1\rangle
=\int_V\left(\pounds_{m_1,p}u\phi_1+h_1(x)|u|^{p-2}u\phi_1-F_u(x,u,v,w)\phi_1\right)d\mu,
$$
$$
\langle\varphi_{w,v}(u,v),\phi_2\rangle
=\int_V\left(\pounds_{m_2,q}v\phi_2+h_2(x)|v|^{q-2}v\phi_2-F_v(x,u,v,w)\phi_2\right)d\mu.
$$
Therefore, the problem of finding the solutions for system (\ref{eq1}) is simplified to seeking the critical points of functional $\varphi_w$ on $W$ (see \cite{Zhang2022}).
\par
Note that $F(x,t,s,w)=\int_0^t F_t(x,t,s,w)dt=\int_0^s F_s(x,t,s,w)ds$ for all $(x,t,s,w)\in V\times \R^2 \times \mathcal{J}$, where $F_t(x,t,s,w)=\frac{\partial F(x,t,s,w)}{\partial t}$ and  $F_s(x,t,s,w)=\frac{\partial F(x,t,s,w)}{\partial s}$. In order to obtain the desired conclusions in the paper, we will assume that $F$ satisfies the following conditions:\\
$(F_1)$ \;  $F(x,0,0,w)=0$ for all $w\in \mathcal{J}$ and $x\in V$, and $F(x,t,s,w)$ is continuously differentiable in $(t,s,w)\in \R^2\times \mathcal{J}$ for all $x\in V$;\\
$(F_2)$ \; $\lim\limits_{|(t,s)|\rightarrow0}\frac{F(x,t,s,w)}{|t|^p+|s|^q}<\min\{\frac{1}{pK_1^p},\frac{1}{qK_2^q}\}$ for all $x\in V$ and $w\in \mathcal{J}$, where
$$
K_1=\frac{\left(\sum_{x\in V}\mu(x)\right)^{\frac{1}{p}}}{\mu_{\min}^\frac{1}{p}h_{1,\min}^\frac{1}{p}},\ \ K_2=\frac{\left(\sum_{x\in V}\mu(x)\right)^{\frac{1}{q}}}{\mu_{\min}^\frac{1}{q}h_{2,\min}^\frac{1}{q}},
$$
where $h_{i,\min}:=\min_{x\in V}h_i(x),i=1,2;$\\
$(F_3)$ \; $\lim\limits_{|(t,s)|\rightarrow\infty}\frac{F(x,t,s,w)}{|t|^p+|s|^q}=+\infty$ for all $x\in V, w\in \mathcal{J}$;\\
$(F_4)$ \;  There exist two constants $\gamma_1>0$ and $\gamma_2>0$ such that
$$
\liminf\limits_{|(t,s)|\rightarrow\infty} \frac{F_t(x,t,s,w)t+F_s(x,t,s,w)s-\max\{p,q\}F(x,t,s,w)}{|t|^{\gamma_1}+|s|^{\gamma_2}}
>0 \; \mbox{ for all} \;\; x\in V,w\in \mathcal{J}.
$$
\vskip2mm
\noindent
{\bf Theorem 3.1.} {\it Assume that $(F_1)$-$(F_4)$ hold. Then system (\ref{eq1}) has at least one mountain pass type solution $(\widetilde{u}_w,\widetilde{v}_w)\not=(0,0)$.}
\vskip2mm
\noindent
{\bf Proof.} When $w$ is fixed, the proof is the same as that in \cite{Zhang2022}. We omit the details.\qed

\vskip2mm
\noindent
\par
In order to obtain the dependence on parameters of solutions, we need the following assumptions:\\
$(H_1)$  There exist a constant $\theta > \max\{p,q\}$  such that
$$\theta F(x,t,s,w)\leq F_t(x,t,s,w)t+F_s(x,t,s,w)s \;\; \text{for all} \;\; x\in V,\ (t,s)\in \R^2\backslash\{(0,0)\},w\in \mathcal{J},$$
$(H_2)$  There exist positive constants $c_1,c_2,r_1$ and $r_2$ with $\min\{r_1,r_2\}>\max\{p,q\}$, such that
$$F_t(x,t,s,w)t+F_s(x,t,s,w)s\leq c_1 |t|^{r_1}+c_2 |s|^{r_2} \;\;\text{for all} \;\; x\in V,(t,s)\in \R^2,w\in \mathcal{J},$$
$(H_3)$  There exist two functions $a\in C(\R,\R^+)$ and $c:V\rightarrow\R^+$ such that
\begin{eqnarray*}
F(x,t,s,w) \geq a(|(t,s)|)c(x) \; \mbox{ for all} \;\; x\in V,(t,s)\in \R^2,w\in \mathcal{J}.
\end{eqnarray*}
\vskip2mm
\noindent
{\bf Remark 3.1.} By $(F_1)$, $(F_2)$, $(H_1)$ and $(H_3)$, we can also obtain that the system (\ref{eq1}) has at least one mountain pass type nontrivial solution. In fact, similar to \cite{Zhang2012},  for any fixed $x\in V, (t,s)\in \R^2\backslash\{(0,0)\}, w\in \mathcal{J}$ and some $M_0>0$, let
$$f(\tau):=f_{x,t,s,w}(\tau)=F(x,\tau t,\tau s,w)\;\;  \text{for all} \;\; \tau\geq \frac{M_0}{|(t,s)|}.$$
$(H_1)$ can be written as
\begin{eqnarray}
\label{eQ18}
\theta F(x,\tau t,\tau s,w)\leq F_{\tau t}(x,\tau t,\tau s,w)\tau t+F_{\tau s}(x,\tau t,\tau s,w)\tau s \;\; \text{for all} \;\; x\in V,(t,s)\in \R^2\backslash\{(0,0)\},w\in \mathcal{J}.
\end{eqnarray}
Thus it follows from the definition of $f(\tau)$ and (\ref{eQ18}) that
\begin{eqnarray*}
f'(\tau)& = & F_{\tau t}(x,\tau t,\tau s,w)t+F_{\tau s}(x,\tau t,\tau s,w)s \\
& = & \frac{1}{\tau}[F_{\tau t}(x,\tau t,\tau s,w)\tau t+F_{\tau s}(x,\tau t,\tau s,w)\tau s]\\
& \geq & \frac{\theta}{\tau}F(x,\tau t,\tau s,w)\\
& =  &  \frac{\theta}{\tau}f(\tau).
\end{eqnarray*}
Let $\rho(\tau)=f'(\tau)-\frac{\theta}{\tau}f(\tau)\geq0$. A simple calculation implies
$$\rho(\tau)=\left[\frac{f(\tau)}{\tau^\theta}\right]'\tau^\theta.$$
Integrating both sides of the above formula, we get
\begin{eqnarray*}
\int^\tau_{\frac{M_0}{|(t,s)|}}\frac{\rho(\tau)}{\tau^\theta}d\tau =\int^\tau_{\frac{M_0}{|(t,s)|}}\left[\frac{f(\tau)}{\tau^\theta}\right]'d\tau = \frac{f(\tau)}{\tau^\theta}-\frac{f\left(\frac{M_0}{|(t,s)|}\right)}{\left(\frac{M_0}{|(t,s)|}\right)^\theta},
\end{eqnarray*}
and then
\begin{eqnarray*}
f(\tau)=\tau^\theta\left(\int^\tau_{\frac{M_0}{|(t,s)|}}\frac{\rho(\tau)}{\tau^\theta}d\tau+\frac{|(t,s)|^\theta f\left(\frac{M_0}{|(t,s)|}\right)}{M_0^\theta}\right),
\end{eqnarray*}
Note that $\int^\tau_{\frac{M_0}{|(t,s)|}}\frac{\rho(\tau)}{\tau^\theta}d\tau \geq 0$. When $|(t,s)|\ge M_0$, we can choose $\tau=1$. Then the above equality implies that
\begin{eqnarray}
\label{eQ19}
f(1)=F(x,t,s,w)\geq \frac{f\left(\frac{M_0}{|(t,s)|}\right)}{M_0^\theta}|(t,s)|^\theta, \ \mbox{for all }|(t,s)|\ge M_0.
\end{eqnarray}
Note that $a\in C(\R,\R^+)$, $c:V\rightarrow\R^+$ and $\theta>\max\{p,q\}$. It follows from (\ref{eQ19}) and $(H_3)$ that
\begin{eqnarray*}
\lim\limits_{|(t,s)|\rightarrow\infty}\frac{F(x,t,s,w)}{|t|^p+|s|^q}&\geq& \lim\limits_{|(t,s)|\rightarrow\infty} \frac{\frac{f\left(\frac{M_0}{|(t,s)|}\right)}{M_0^\theta}|(t,s)|^\theta}{|t|^p+|s|^q}=\lim\limits_{|(t,s)|\rightarrow\infty} \frac{\frac{F\left(x,\frac{M_0}{|(t,s)|}t,\frac{M_0}{|(t,s)|}s,w\right)}{M_0^\theta}|(t,s)|^\theta}{|t|^p+|s|^q} \\
&\geq &\lim\limits_{|(t,s)|\rightarrow\infty} \frac{a(|(\frac{M_0}{|(t,s)|}t,\frac{M_0}{|(t,s)|}s)|)c(x)}{M_0^\theta}\frac{|(t,s)|^\theta}{|t|^p+|s|^q}\\
&\geq & \min\limits_{|(t,s)|=M_0}a(|(t,s)|)\min\limits_{x\in V}c(x)\frac{1}{M_0^\theta}\lim\limits_{|(t,s)|\rightarrow\infty} \frac{|(t,s)|^\theta}{|t|^p+|s|^q}\\
&=    &+\infty
\end{eqnarray*}
for all $x\in V,w\in \mathcal{J}$. So $(F_3)$ holds.
Next, by $(H_1)$, we can obtain
\begin{eqnarray*}
F_t(x,t,s,w)t+F_s(x,t,s,w)s \geq \theta F(x,t,s,w)=(\theta-\max\{p,q\})F(x,t,s,w)+\max\{p,q\}F(x,t,s,w),
\end{eqnarray*}
for all  $x\in V, (t,s)\in \R^2\backslash\{(0,0)\}, w\in \mathcal{J}$. Let $0<\gamma_1\leq p$ and $0<\gamma_2\leq q$. Then, by $(F_3)$, we infer that
\begin{eqnarray*}
&&\liminf\limits_{|(t,s)|\rightarrow\infty} \frac{F_t(x,t,s,w)t+F_s(x,t,s,w)s-\max\{p,q\}F(x,t,s,w)}{|t|^{\gamma_1}+|s|^{\gamma_2}}\\
& \geq &\liminf\limits_{|(t,s)|\rightarrow\infty}\frac{(\theta-\max\{p,q\})F(x,t,s,w)}{|t|^{\gamma_1}+|s|^{\gamma_2}}\\
& \geq & (\theta-\max\{p,q\})\liminf\limits_{|(t,s)|\rightarrow\infty}\frac{F(x,t,s,w)}{|t|^p+|s|^q}\\
& =    & +\infty
\end{eqnarray*}
for all $x\in V,w\in \mathcal{J}$.  Hence, $(F_4)$ holds.\qed

\vskip2mm
\noindent
{\bf Theorem 3.2.} {\it Assume that $(F_1)$, $(F_2)$ and $(H_1)$-$(H_3)$ hold. Then, for any $w\in \mathcal{J}$, system (\ref{eq1}) has at least one nontrivial solution $(u_w,v_w)$ satisfying
\begin{eqnarray}
\label{EQ4}
  C_1 \leq \|(u_w,v_w)\| \leq C_2,
\end{eqnarray}
where $C_1$ and $C_2$ are positive constants which are independent on $w$ and the values are given in $(\ref{EQ17})$ and $(\ref{EQ10})$ below,  respectively. }
\vskip2mm
\noindent
{\bf Proof.} Choose arbitrary but fixed $w\in \mathcal{J}$. According to Theorem 3.1 and Remark 3.1, it is obvious that system (\ref{eq1}) has at least one nontrivial solution $(u_w,v_w)$.  Furthermore, from the proof of Lemma 3.3 in  \cite{Zhang2022}, there exists a $(0,0)\neq(u_{0},v_{0})\in W$  which is independent on $w$ since $(F_3)$ holds for all $w\in \mathcal{J}$,  such that
\begin{eqnarray}
\label{EQ5}
  \varphi_w(u_w,v_w)=\inf_{\gamma\in\Gamma} \max_{t\in[0,1]} \varphi_w(\gamma(t)):=C^*>0,
\end{eqnarray}
where
$$
\Gamma:=\{\gamma\in C([0,1],W):\gamma(0)=(0,0),\gamma(1)=(u_{0},v_{0})\}.
$$
Since $(u_w,v_w)$ is a weak solution of (\ref{eq1}), then
\begin{small}
\begin{eqnarray}
\label{EQ6}
\|u_w\|^p_{W^{m_1,p}(V)}+\|v_w\|^p_{W^{m_2,q}(V)}-\int_ VF_{u_w}(x,u_w(x),v_w(x),w)u_w(x) d\mu
-\int_ VF_{v_w}(x,u_w(x),v_w(x),w)v_w(x)d\mu=0.
\end{eqnarray}
\end{small}Next, we first prove the first inequality in (\ref{EQ4}). Without loss of generality, we assume that $p\geq q$. We divide the proof into the following four cases.\\
$(1)$ Assume that $\|u_w\|_{W^{m_1,p}(V)} >1$ and $\|v_w\|_{W^{m_2,q}(V)}>1$. By  $(H_2)$, (\ref{EQ6}) and Lemma 2.1, we have
\begin{eqnarray}
\|(u_w,v_w)\|^q
\notag
&  =  & (\|u_w\|_{W^{m_1,p}(V)}+\|v_w\|_{W^{m_2,q}(V)})^q\\
\notag
& \leq & 2^{q-1}(\|u_w\|^q_{W^{m_1,p}(V)}+\|v_w\|^q_{W^{m_2,q}(V)}) \\
\notag
& \leq & 2^{q-1}(\|u_w\|^p_{W^{m_1,p}(V)}+\|v_w\|^q_{W^{m_2,q}(V)}) \\
\notag
&  =  & 2^{q-1}\left(\int_ VF_{u_w}(x,u_w(x),v_w(x),w)u_w(x) d\mu+\int_ VF_{v_w}(x,u_w(x),v_w(x),w)v_w(x)d\mu \right)\\
\notag
& \leq & 2^{q-1}\left(\int_ V c_1 |u_w(x)|^{r_1} d\mu+\int_ V c_2 |v_w(x)|^{r_2}d\mu \right)\\
\notag
& \leq & 2^{q-1} |V|\max{\{c_1,c_2\}}\max{\{b^{r_1},d^{r_2}\}}(\|u_w\|_{W^{m_1,p}(V)}^{r_1}+ \|v_w\|_{W^{m_2,q}(V)}^{r_2})\\
\notag
& \leq & 2^{q-1} |V|\max{\{c_1,c_2\}}\max{\{b^{r_1},d^{r_2}\}}\|(u_w,v_w)\|^{\max\{r_1,r_2\}},
\end{eqnarray}
which implies that
\begin{eqnarray}
\|(u_w,v_w)\|\geq \left({\frac{1}{2^{q-1} |V|\max{\{c_1,c_2\}}\max{\{b^{r_1},d^{r_2}\}}}}\right)^{\frac{1}{\max\{r_1,r_2\}-q}}.
\notag
\end{eqnarray}
$(2)$ Assume that $\|u_w\|_{W^{m_1,p}(V)}\leq1$ and $\|v_w\|_{W^{m_2,q}(V)}\leq1$. According to $(H_2)$, (\ref{EQ6}) and Lemma 2.1, there exists
\begin{eqnarray}
\|(u_w,v_w)\|^p
\notag
&  =  & (\|u_w\|_{W^{m_1,p}(V)}+\|v_w\|_{W^{m_2,q}(V)})^p\\
\notag
& \leq & 2^{p-1}(\|u_w\|^p_{W^{m_1,p}(V)}+\|v_w\|^p_{W^{m_2,q}(V)}) \\
\notag
& \leq & 2^{p-1}(\|u_w\|^p_{W^{m_1,p}(V)}+\|v_w\|^q_{W^{m_2,q}(V)}) \\
\notag
&  =  & 2^{p-1}\left(\int_ VF_{u_w}(x,u_w(x),v_w(x),w)u_w(x) d\mu+\int_ VF_{v_w}(x,u_w(x),v_w(x),w)v_w(x)d\mu \right)\\
\notag
& \leq & 2^{p-1}\left(\int_ V c_1 |u_w(x)|^{r_1} d\mu+\int_ V c_2 |v_w(x)|^{r_2}d\mu \right)\\
\notag
& \leq & 2^{p-1} |V|\max{\{c_1,c_2\}}\max{\{b^{r_1},d^{r_2}\}}(\|u_w\|_{W^{m_1,p}(V)}^{r_1}+ \|v_w\|_{W^{m_2,q}(V)}^{r_2})\\
\notag
& \leq & 2^{p-1} |V|\max{\{c_1,c_2\}}\max{\{b^{r_1},d^{r_2}\}}\|(u_w,v_w)\|^{\min\{r_1,r_2\}},
\end{eqnarray}
which implies that
\begin{eqnarray}
\|(u_w,v_w)\|\geq \left({\frac{1}{2^{p-1} |V|\max{\{c_1,c_2\}}\max{\{b^{r_1},d^{r_2}\}}}}\right)^{\frac{1}{\min\{r_1,r_2\}-p}}.
\notag
\end{eqnarray}
$(3)$ Assume that $\|u_w\|_{W^{m_1,p}(V)}>1$ and $\|v_w\|_{W^{m_2,q}(V)}\leq 1$. Then $\|(u_w,v_w)\|>1$, and thus by Lemma 2.1,
\begin{eqnarray}
\|(u_w,v_w)\|^q
\notag
&  =  & (\|u_w\|_{W^{m_1,p}(V)}+\|v_w\|_{W^{m_2,q}(V)})^q\\
\notag
& \leq & 2^{q-1}(\|u_w\|^q_{W^{m_1,p}(V)}+\|v_w\|^q_{W^{m_2,q}(V)}) \\
\notag
& \leq & 2^{q-1}(\|u_w\|^p_{W^{m_1,p}(V)}+\|v_w\|^q_{W^{m_2,q}(V)}) \\
\notag
&  =  & 2^{q-1}\left(\int_ VF_{u_w}(x,u_w(x),v_w(x),w)u_w(x) d\mu+\int_ VF_{v_w}(x,u_w(x),v_w(x),w)v_w(x)d\mu \right)\\
\notag
& \leq & 2^{q-1}\left(\int_ V c_1 |u_w(x)|^{r_1} d\mu+\int_ V c_2 |v_w(x)|^{r_2}d\mu \right)\\
\notag
& \leq & 2^{q-1} |V|\max{\{c_1,c_2\}}\max{\{b^{r_1},d^{r_2}\}}(\|u_w\|_{W^{m_1,p}(V)}^{r_1}+ \|v_w\|_{W^{m_2,q}(V)}^{r_2})\\
\notag
& \leq & 2^{q-1} |V|\max{\{c_1,c_2\}}\max{\{b^{r_1},d^{r_2}\}}(\|(u_w,v_w)\|^{r_1}+ \|(u_w,v_w)\|^{r_2})\\
\notag
& \leq & 2^{q-1} |V|\max{\{c_1,c_2\}}\max{\{b^{r_1},d^{r_2}\}}2\|(u_w,v_w)\|^{\max\{r_1,r_2\}}.
\end{eqnarray}
Then we have
\begin{eqnarray}
\|(u_w,v_w)\|\geq \left({\frac{1}{2^q |V|\max{\{c_1,c_2\}}\max{\{b^{r_1},d^{r_2}\}}}}\right)^{\frac{1}{\max\{r_1,r_2\}-q}}.
\notag
\end{eqnarray}
$(4)$ Assume that $\|u_w\|_{W^{m_1,p}(V)}\leq1$ and $\|v_w\|_{W^{m_2,q}(V)}>1$. Note that
\begin{eqnarray}
\label{EQ7}
\|u_w\|^p_{W^{m_1,p}(V)}+\|v_w\|^q_{W^{m_2,q}(V)}
\notag
& \geq & \|v_w\|_{W^{m_2,q}(V)}^{q-p}(\|u_w\|^p_{W^{m_1,p}(V)}+\|v_w\|^p_{W^{m_2,q}(V)})\\
\notag
& \geq & \|(u_w,v_w)\|^{q-p}(\|u_w\|^p_{W^{m_1,p}(V)}+\|v_w\|^p_{W^{m_2,q}(V)})\\
& \geq & \frac{\|(u_w,v_w)\|^{q-p}}{2^{p-1}}\|(u_w,v_w)\|^p = \frac{\|(u_w,v_w)\|^q}{2^{p-1}}.
\end{eqnarray}
On the other hand, it follows from Lemma 2.1 that
\begin{eqnarray}
\|u_w\|^p_{W^{m_1,p}(V)}+\|v_w\|^q_{W^{m_2,q}(V)}
\notag
&  =  & \int_ VF_{u_w}(x,u_w(x),v_w(x),w)u_w(x) d\mu+\int_ VF_{v_w}(x,u_w(x),v_w(x),w)v_w(x)d\mu\\
\notag
& \leq & \int_ V c_1 |u_w(x)|^{r_1} d\mu+\int_ V c_2 |v_w(x)|^{r_2}d\mu\\
\notag
& \leq & |V|\max{\{c_1,c_2\}}\max{\{b^{r_1},d^{r_2}\}}(\|u_w\|_{W^{m_1,p}(V)}^{r_1}+ \|v_w\|_{W^{m_2,q}(V)}^{r_2})\\
\notag
& \leq & |V|\max{\{c_1,c_2\}}\max{\{b^{r_1},d^{r_2}\}}(\|(u_w,v_w)\|^{r_1}+ \|(u_w,v_w)\|^{r_2})\\
\notag
& \leq & |V|\max{\{c_1,c_2\}}\max{\{b^{r_1},d^{r_2}\}}2\|(u_w,v_w)\|^{\max\{r_1,r_2\}}.
\end{eqnarray}
Then, combining with (\ref{EQ7}), we can derive that
\begin{eqnarray}
\|(u_w,v_w)\| \geq \left({\frac{1}{2^{p} |V|\max{\{c_1,c_2\}}\max{\{b^{r_1},d^{r_2}\}}}}\right)^{\frac{1}{\max\{r_1,r_2\}-q}}.
\notag
\end{eqnarray}
Let
\begin{small}
\begin{eqnarray*}
\label{EQ11}
A_1:=
\min\left\{\left({\frac{1}{2^p |V|\max{\{c_1,c_2\}}\max{\{b^{r_1},d^{r_2}\}}}}\right)^{\frac{1}{\max\{r_1,r_2\}-q}},
\left({\frac{1}{2^{p-1} |V|\max{\{c_1,c_2\}}\max{\{b^{r_1},d^{r_2}\}}}}\right)^{\frac{1}{\min\{r_1,r_2\}-p}}\right\}.
\end{eqnarray*}
\end{small}Then, $\|(u_w,v_w)\| \geq A_1$.\\
For the case of $q>p$, we can also obtain $\|(u_w,v_w)\| \geq A_2$, where
\begin{small}
\begin{eqnarray*}
\label{EQ8}
\notag
A_2:=\min\left\{\left({\frac{1}{2^q |V|\max{\{c_1,c_2\}}\max{\{b^{r_1},d^{r_2}\}}}}\right)^{\frac{1}{\max\{r_1,r_2\}-p}},\right.\\
\left.\left({\frac{1}{2^{q-1} |V|\max{\{c_1,c_2\}}\max{\{b^{r_1},d^{r_2}\}}}}\right)^{\frac{1}{\min\{r_1,r_2\}-q}}\right\}.
\end{eqnarray*}
\end{small}
Furthermore, we have
\begin{eqnarray}
\label{EQ17}
\|(u_w,v_w)\| \geq C_1:=\min\{A_1,A_2\}.
\end{eqnarray}
Next, we prove the second inequality in (\ref{EQ4}). Note that
\begin{eqnarray}
\notag
 \varphi_w(u_w,v_w)=\frac{1}{p} \|u_w\|_{W^{m_1,p}(V)}^p+\frac{1}{q} \|v_w\|_{W^{m_2,q}(V)}^q-\int_V F(x,u_w,v_w,w)d\mu.
\end{eqnarray}
Then we can get
\begin{eqnarray}
\label{EQ13}
\theta \varphi_w(u_w,v_w)+\theta \int_V F(x,u_w(x),v_w(x),w)d\mu=\frac{\theta}{p}\|u_w\|_{W^{m_1,p}(V)}^p+\frac{\theta}{q} \|v_w\|_{W^{m_2,q}(V)}^q.
\end{eqnarray}
It follows from (\ref{EQ6}) and (\ref{EQ13}) that
\begin{eqnarray*}
& &    \left(\frac{\theta}{p}-1\right)\|u_w\|_{W^{m_1,p}(V)}^p+\left(\frac{\theta}{q}-1\right)\|v_w\|_{W^{m_2,q}(V)}^q\\
&=&\theta \varphi_w(u_w,v_w)+\theta \int_V F(x,u_w(x),v_w(x),w)d\mu \\
&&-\int_VF_{u_w}(x,u_w(x),v_w(x),w)u_w(x)d\mu -\int_VF_{v_w}(x,u_w(x),v_w(x),w)v_w(x)d\mu.
\end{eqnarray*}
From $(H_1)$, we can get that
\begin{eqnarray}
\label{EQ9}
\left(\frac{\theta}{p}-1\right)\|u_w\|_{W^{m_1,p}(V)}^p+\left(\frac{\theta}{q}-1\right)\|v_w\|_{W^{m_2,q}(V)}^q \leq \theta \varphi_w(u_w,v_w).
\end{eqnarray}
We need divide the proof into the following four cases, which is referred to \cite{Pang2023}. Without loss of generality, we assume that $p\geq q$.\\
$(1)$ Assume that $\|u_w\|_{W^{m_1,p}(V)}>1$ and $\|v_w\|_{W^{m_2,q}(V)}>1$. Then we have
\begin{eqnarray}
\left(\frac{\theta}{p}-1\right)\|u_w\|^p_{W^{m_1,p}(V)}+\left(\frac{\theta}{q}-1\right)\|v_w\|^q_{W^{m_2,q}(V)}
\notag
& \geq & \left(\frac{\theta}{p}-1\right)\left(\|u_w\|^q_{W^{m_1,p}(V)}+\|v_w\|^q_{W^{m_2,q}(V)}\right) \\
\notag
& \geq & \frac{(\frac{\theta}{p}-1)}{2^{q-1}}\|(u_w,v_w)\|^q\\
\notag
& \geq & \frac{(\frac{\theta}{p}-1)}{2^{p-1}}\|(u_w,v_w)\|^q,
\end{eqnarray}
together with (\ref{EQ9}), which implies that
\begin{eqnarray}
\frac{(\frac{\theta}{p}-1)}{2^{p-1}}\|(u_w,v_w)\|^q \leq \theta \varphi_w(u_w,v_w).
\notag
\end{eqnarray}
$(2)$ Assume that $\|u_w\|_{W^{m_1,p}(V)}\leq1$ and $\|v_w\|_{W^{m_2,q}(V)}\leq1$. Then we have
\begin{eqnarray}
\left(\frac{\theta}{p}-1\right)\|u_w\|^p_{W^{m_1,p}(V)}+\left(\frac{\theta}{q}-1\right)\|v_w\|^q_{W^{m_2,q}(V)}
\notag
& \geq & \left(\frac{\theta}{p}-1\right)\left(\|u_w\|^p_{W^{m_1,p}(V)}+\|v_w\|^p_{W^{m_2,q}(V)}\right) \\
\notag
& \geq & \frac{(\frac{\theta}{p}-1)}{2^{p-1}}\|(u_w,v_w)\|^p\\
\notag
& \geq &
\begin{cases}
\frac{(\frac{\theta}{p}-1)}{2^{p-1}}\|(u_w,v_w)\|^p,\;\;\;\|(u_w,v_w)\|\leq1, \\
\frac{(\frac{\theta}{p}-1)}{2^{p-1}}\|(u_w,v_w)\|^q,\;\;\;1<\|(u_w,v_w)\|\leq2.
\end{cases}
\end{eqnarray}
Then by (\ref{EQ9}),
\begin{eqnarray}
\begin{cases}
\frac{(\frac{\theta}{p}-1)}{2^{p-1}}\|(u_w,v_w)\|^p \leq \theta \varphi_w(u_w,v_w),\;\;\;\|(u_w,v_w)\|\leq1,\\
\notag
\frac{(\frac{\theta}{p}-1)}{2^{p-1}}\|(u_w,v_w)\|^q \leq \theta \varphi_w(u_w,v_w),\;\;\;1<\|(u_w,v_w)\|\leq2.
\notag
\end{cases}
\end{eqnarray}
$(3)$ Assume that $\|u_w\|_{W^{m_1,p}(V)}>1$ and $\|v_w\|_{W^{m_2,q}(V)}\leq1$. Then we have
\begin{eqnarray}
\left(\frac{\theta}{p}-1\right)\|u_w\|^p_{W^{m_1,p}(V)}+\left(\frac{\theta}{q}-1\right)\|v_w\|^q_{W^{m_2,q}(V)}
\notag
& \geq & \left(\frac{\theta}{p}-1\right)\left(\|u_w\|^q_{W^{m_1,p}(V)}+\|v_w\|^q_{W^{m_2,q}(V)}\right) \\
\notag
& \geq & \frac{(\frac{\theta}{p}-1)}{2^{q-1}}\|(u_w,v_w)\|^q\\
\notag
& \geq & \frac{(\frac{\theta}{p}-1)}{2^{p-1}}\|(u_w,v_w)\|^q.
\end{eqnarray}
Furthermore, according to (\ref{EQ9}),
\begin{eqnarray}
\frac{(\frac{\theta}{p}-1)}{2^{p-1}}\|(u_w,v_w)\|^q \leq \theta \varphi_w(u_w,v_w).
\notag
\end{eqnarray}
$(4)$ Assume that $\|u_w\|_{W^{m_1,p}(V)}\leq1$ and $\|v_w\|_{W^{m_2,q}(V)}>1$. Then we have
\begin{eqnarray}
\notag
&&\left(\frac{\theta}{p}-1\right)\|u_w\|^p_{W^{m_1,p}(V)}+\left(\frac{\theta}{q}-1\right)\|v_w\|^q_{W^{m_2,q}(V)}\\
\notag
& \geq & \min\left\{\left(\frac{\theta}{p}-1\right),\left(\frac{\theta}{q}-1\right)\|v_w\|^{q-p}_{W^{m_2,q}(V)}\right\}\left(\|u_w\|^p_{W^{m_1,p}(V)}+\|v_w\|^p_{W^{m_2,q}(V)}\right) \\
\notag
& \geq &
\min\left\{\left(\frac{\theta}{p}-1\right),\left(\frac{\theta}{q}-1\right)\|(u_w,v_w)\|^{q-p}\right\}\left(\|u_w\|^p_{W^{m_1,p}(V)}+\|v_w\|^p_{W^{m_2,q}(V)}\right) \\
\notag
& \geq &
\min\left\{\left(\frac{\theta}{p}-1\right),\left(\frac{\theta}{q}-1\right)\|(u_w,v_w)\|^{q-p}\right\}\frac{1}{2^{p-1}}\|(u_w,v_w)\|^p \\
\notag
& = &
\min\left\{\frac{(\frac{\theta}{p}-1)}{2^{p-1}}\|(u_w,v_w)\|^p,\frac{(\frac{\theta}{q}-1)}{2^{p-1}}\|(u_w,v_w)\|^q\right\} \\
\notag
& \geq & \min\left\{\frac{(\frac{\theta}{p}-1)}{2^{p-1}}\|(u_w,v_w)\|^q,\frac{(\frac{\theta}{q}-1)}{2^{p-1}}\|(u_w,v_w)\|^q\right\} \\
& = & \frac{(\frac{\theta}{p}-1)}{2^{p-1}}\|(u_w,v_w)\|^q.
\end{eqnarray}
Hence,
\begin{eqnarray}
\notag
\frac{(\frac{\theta}{p}-1)}{2^{p-1}}\|(u_w,v_w)\|^q \leq \theta \varphi_w(u_w,v_w).
\notag
\end{eqnarray}
From the above formula we derive
\begin{eqnarray}
\begin{cases}
\frac{(\frac{\theta}{p}-1)}{2^{p-1}}\|(u_w,v_w)\|^{p} \leq \theta \varphi_w(u_w,v_w),\;\;\;\|(u_w,v_w)\|\leq1,\\
\notag
\frac{(\frac{\theta}{p}-1)}{2^{p-1}}\|(u_w,v_w)\|^{q} \leq \theta \varphi_w(u_w,v_w),\;\;\;\|(u_w,v_w)\|>1.
\notag
\end{cases}
\end{eqnarray}
By (\ref{EQ5}), $(H_3)$ and the definition of $\varphi_\omega$,
\begin{eqnarray}\label{aa2}
\notag
\theta \varphi_w(u_w,v_w)
&  =  & \theta \inf_{\gamma\in\Gamma} \max_{t\in[0,1]} \varphi_w(\gamma(t)) \\
\notag
& \leq & \theta \max_{t\in[0,1]} \varphi_w (t u_0,t v_0) \\
\notag
& = & \theta \max_{t\in[0,1]} \left(\frac{1}{p}\|t u_0\|^p_{W^{m_1,p}(V)}+\frac{1}{q}\|t v_0\|^q_{W^{m_2,q}(V)} - \int_V F(x,t u_0, t v_0,w)d\mu\right) \\
\notag
& \leq & \theta \max_{t\in[0,1]} \left(\frac{1}{p}t^p\|u_0\|^p_{W^{m_1,p}(V)}+\frac{1}{q}t^q\| v_0\|^q_{W^{m_2,q}(V)}\right)\\
& \leq &  \theta  \left(\frac{1}{p}\|u_0\|^p_{W^{m_1,p}(V)}+\frac{1}{q}\| v_0\|^q_{W^{m_2,q}(V)}\right).
\end{eqnarray}
Then,
\begin{eqnarray}
\begin{cases}
\label{eq8}
\frac{(\frac{\theta}{p}-1)}{2^{p-1}}\|(u_w,v_w)\|^{p} \leq \theta \varphi_w(u_w,v_w)\leq \theta \left(\frac{1}{p}\|u_0\|^p_{W^{m_1,p}(V)}+\frac{1}{q}\| v_0\|^q_{W^{m_2,q}(V)}\right) ,\;\;\;\|(u_w,v_w)\|\leq1,\\
\notag
\frac{(\frac{\theta}{p}-1)}{2^{p-1}}\|(u_w,v_w)\|^{q} \leq \theta \varphi_w(u_w,v_w)\leq \theta \left(\frac{1}{p}\|u_0\|^p_{W^{m_1,p}(V)}+\frac{1}{q}\| v_0\|^q_{W^{m_2,q}(V)}\right),\;\;\;\|(u_w,v_w)\|>1.
\notag
\end{cases}
\end{eqnarray}
Therefore,
\begin{eqnarray*}
\label{EQ14}
\begin{cases}
\label{eq8}
\|(u_w,v_w)\|\leq  \left({\frac{p\theta 2^{p-1}\left(\frac{1}{p}\|u_0\|^p_{W^{m_1,p}(V)}+\frac{1}{q}\| v_0\|^q_{W^{m_2,q}(V)}\right)}{\theta-p}}\right)^{\frac{1}{p}}:=A_3,\;\;\;\;\ \|(u_w,v_w)\|\leq1,\\
\|(u_w,v_w)\|\leq  \left({\frac{p\theta 2^{p-1}\left(\frac{1}{p}\|u_0\|^p_{W^{m_1,p}(V)}+\frac{1}{q}\| v_0\|^q_{W^{m_2,q}(V)}\right)}{\theta-p}}\right)^{\frac{1}{q}}:=A_4,\;\;\;\;\ \|(u_w,v_w)\|>1.
\end{cases}
\end{eqnarray*}
Similarly, with $q>p$, we can also get,
\begin{eqnarray*}
\label{EQ14}
\begin{cases}
\label{eq8}
\|(u_w,v_w)\|\leq  \left({\frac{q\theta 2^{q-1}\left(\frac{1}{p}\|u_0\|^p_{W^{m_1,p}(V)}+\frac{1}{q}\| v_0\|^q_{W^{m_2,q}(V)}\right)}{\theta-q}}\right)^{\frac{1}{q}}:=A_5,\;\;\;\;\ \|(u_w,v_w)\|\leq1,\\
\|(u_w,v_w)\|\leq  \left({\frac{q\theta 2^{q-1}\left(\frac{1}{p}\|u_0\|^p_{W^{m_1,p}(V)}+\frac{1}{q}\| v_0\|^q_{W^{m_2,q}(V)}\right)}{\theta-q}}\right)^{\frac{1}{p}}:=A_6,\;\;\;\;\ \|(u_w,v_w)\|>1.
\end{cases}
\end{eqnarray*}
So,
\begin{eqnarray}
\label{EQ10}
\|(u_w,v_w)\|\leq C_2:=\max\{A_3,A_4,A_5,A_6\}.
\end{eqnarray}
In conclusion, (\ref{EQ4}) holds and the proof is completed.\qed
\vskip2mm
\par
The  result concerning with the dependence on the parameters of solutions is stated as follows.
\vskip2mm
\noindent
{\bf Theorem 3.3.} {\it Assume that  $(F_1)$, $(F_2)$ and $(H_1)$-$(H_3)$ are satisfied. Let $\{w_n\} \subset \mathcal{J}$ be a convergent sequence of parameters with $\lim\limits_{n\rightarrow\infty}w_n=w_0$. For any nontrivial solutions sequence $\{(u_n,v_n)\}$ related to $\{w_n\}$ in (\ref{eq1}), there exists a subsequence $\{(u_{n_i},v_{n_i})\}\subset W$ such that $\lim\limits_{i\rightarrow\infty}(u_{n_i},v_{n_i})=(u_0,v_0)\in W$. Moreover, $(u_0,v_0)$ is a nontrivial solution to (\ref{eq1}) corresponding to $w_0$, and $C_1\leq\|(u_0,v_0)\|\leq C_2$.}
\vskip2mm
\noindent
{\bf Proof.} For each $n\in \mathbb N$, we take $(u_n,v_n)\in W$ as a solution to (\ref{eq1}) corresponding $w=w_n$. By Theorem 3.2, we have $C_1\leq\|(u_n,v_n)\|\leq C_2$ for $n=1,2,...$, where $C_1,C_2$ are given by (\ref{EQ17}) and (\ref{EQ10}). Then by virtue of the fact that $W$ is a finite dimensional space, there is a subsequence of $\{(u_n,v_n)\}$, still denoted by $\{(u_n,v_n)\}$, and a $(u_0,v_0)\in W$ such that $(u_n,v_n)\rightarrow (u_0,v_0)$ in $W$ as $n\to \infty$.
\par
By Lemma 2.1, we get
\begin{eqnarray}
\label{EQ16}
u_n(x)\rightarrow u_0(x), v_n(x)\rightarrow v_0(x) \;\;\;\;\text{for all}\;\;x\in V.
\end{eqnarray}
Since $(u_n,v_n)$ is a solution to (\ref{eq1}), then for any $(\phi_1,\phi_2) \in W$, it holds that
\begin{eqnarray}
\label{EQ1}
 \int_V\left(\pounds_{m_1,p}u_n\phi_1+h_1(x)|u_n|^{p-2}u_n\phi_1-F_u(x,u_n,v_n,w_n)\phi_1\right)d\mu\nonumber\\
 +\int_V\left(\pounds_{m_2,q}v_n\phi_2+h_2(x)|v_n|^{q-2}v_n\phi_2-F_v(x,u_n,v_n,w_n)\phi_2\right)d\mu=0.
\end{eqnarray}
\par
Since $F$ is continuous, by (\ref{EQ16}), we obtain $F_u (\cdot,u_n(\cdot),v_n(\cdot),w_n)\rightarrow F_u(\cdot,u_0(\cdot),v_0(\cdot),w_0)$ and $F_v (\cdot,u_n(\cdot),v_n(\cdot),w_n)\\
\rightarrow F_v(\cdot,u_0(\cdot),v_0(\cdot),w_0)$ in $W$.
By the definition of $\pounds_{m,p}$, (\ref{Eq1})-(\ref{Eq4}), the fact that the sum of finite terms is still limited, and (\ref{EQ16}),  when $m_1=1$, as $n \rightarrow \infty$, we have
\begin{eqnarray*}
|\nabla u_n(x)|=\sqrt{\Gamma(u_n)(x)}
 = \left(\frac{1}{2\mu(x)}\sum\limits_{y\thicksim x}w_{xy}(u_n(y)-u_n(x))^2\right)^{\frac{1}{2}}
 \rightarrow  \left(\frac{1}{2\mu(x)}\sum\limits_{y\thicksim x}w_{xy}(u_0(y)-u_0(x))^2\right)^{\frac{1}{2}}=|\nabla u_0(x)|.
\end{eqnarray*}
When $m_1=2$, we have
\begin{eqnarray*}
|\nabla^2 u_n(x)|=|\Delta u_n(x)|=\left|\frac{1}{\mu(x)}\sum\limits_{y\thicksim x}w_{xy}(u_n(y)-v_n(x))\right|\rightarrow \left|\frac{1}{\mu(x)}\sum\limits_{y\thicksim x}w_{xy}(u_0(y)-v_0(x))\right|=|\Delta u_0(x)|=|\nabla^2 u_0(x)|.
\end{eqnarray*}
When $m_1=3$, we have
\begin{eqnarray*}
|\nabla^3 u_n(x)|=|\nabla \Delta u_n(x)|=\sqrt{\Gamma(\Delta u_n)(x)}
& = &\left(\frac{1}{2\mu(x)}\sum\limits_{y\thicksim x}w_{xy}(\Delta u_n(y)-\Delta u_n(x))^2\right)^{\frac{1}{2}} \\
& \rightarrow & \left(\frac{1}{2\mu(x)}\sum\limits_{y\thicksim x}w_{xy}(\Delta u_0(y)-\Delta u_0(x))^2\right)^{\frac{1}{2}}=|\nabla^3 u_0(x)|.
\end{eqnarray*}
When $m_1=4$, we have
\begin{eqnarray*}
|\nabla^4 u_n(x)|=|\Delta^2 u_n(x)|=|\Delta(\Delta u_n(x)|
& = & \left|\frac{1}{\mu(x)}\sum\limits_{y\thicksim x}w_{xy}(\Delta u_n(y)-\Delta v_n(x))\right|\\
& \rightarrow & \left|\frac{1}{\mu(x)}\sum\limits_{y\thicksim x}w_{xy}(\Delta u_0(y)-\Delta v_0(x))\right|=|\nabla^4 u_0(x)|.
\end{eqnarray*}
Furthermore, the mathematical induction implies that as $n\rightarrow \infty$,
\begin{eqnarray*}
 \begin{cases}
|\nabla\Delta^{\frac{m_1-1}{2}}u_n| \rightarrow |\nabla\Delta^{\frac{m_1-1}{2}}u_0|,& \;\;\;\text {when $m$ is odd},\\
|\Delta^{\frac{m_1}{2}}u_n| \rightarrow |\Delta^{\frac{m_1}{2}}u_0|,&\;\;\;\text {when $m$ is even}.
 \end{cases}
\end{eqnarray*}
So when $m$ is odd, we can obtain
\begin{eqnarray*}
&&\int_V\pounds_{m_1,p}u_n\phi_1d\mu \\
& = &  \int_V |\nabla \Delta^{\frac{m_1-1}{2}}u_n|^{p-2} \frac{1}{2\mu(x)}
      \sum\limits_{y\thicksim x}w_{xy}\left(\Delta^{\frac{m_1-1}{2}}u_n(y)-\Delta^{\frac{m_1-1}{2}}u_n(x)\right)
      \left(\Delta^{\frac{m-1}{2}}\phi_1(y)-\Delta^{\frac{m-1}{2}}\phi_1(x)\right)d\mu\\
& \rightarrow &  \int_V |\nabla \Delta^{\frac{m_1-1}{2}}u_0|^{p-2} \frac{1}{2\mu(x)}
      \sum\limits_{y\thicksim x}w_{xy}\left(\Delta^{\frac{m_1-1}{2}}u_0(y)-\Delta^{\frac{m_1-1}{2}}u_0(x)\right)
     \left(\Delta^{\frac{m-1}{2}}\phi_1(y)-\Delta^{\frac{m-1}{2}}\phi_1(x)\right)d\mu\\
& = & \int_V\pounds_{m_1,p}u_0\phi_1d\mu.
\end{eqnarray*}
When $m$ is even, there holds
\begin{eqnarray}
\int_V\pounds_{m_1,p}u_n\phi_1d\mu
\notag
& = &  \int_V|\Delta^{\frac{m_1}{2}}u_n|^{p-2}\Delta^{\frac{m_1}{2}}u_n\Delta^{\frac{m_1}{2}}\phi_1 d\mu \\
& \rightarrow &  \int_V|\Delta^{\frac{m_1}{2}}u_0|^{p-2}\Delta^{\frac{m_1}{2}}u_0\Delta^{\frac{m_1}{2}}\phi_1 d\mu
=\int_V\pounds_{m_1,p}u_0\phi_1d\mu.
\end{eqnarray}
In a word, from the above, we can get that $\int_V\pounds_{m_1,p}u_n\phi_1d\mu\rightarrow \int_V\pounds_{m_1,p}u_0\phi_1d\mu$. Similarly, we also have $\int_V\pounds_{m_2,q}v_n\phi_2d\mu\\
\rightarrow\int_V\pounds_{m_2,q}v_0\phi_2d\mu$. Moreover, by (\ref{EQ16}), as $n \rightarrow \infty$, we also get
\begin{eqnarray*}
\begin{cases}
  \int_V h_1(x)|u_n|^{p-2}u_n\phi_1d\mu \rightarrow \int_Vh_1(x)|u_0|^{p-2}u_0\phi_1d\mu,\\
  \int_V h_2(x)|v_n|^{q-2}v_n\phi_2d\mu \rightarrow \int_Vh_2(x)|v_0|^{q-2}v_0\phi_2d\mu,
\end{cases}
\text{for any}\;\;(\phi_1,\phi_2) \in W.
\end{eqnarray*}
Thus combining with (\ref{EQ1}), we obtian that for any $(\phi_1,\phi_2) \in W$,
\begin{small}
\begin{eqnarray*}
  \int_V\left(\pounds_{m_1,p}u_0\phi_1+h_1(x)|u_0|^{p-2}u_0\phi_1-F_{u_0}(x,u_0,v_0,w_0)\phi_1\right)d\mu
 +\int_V\left(\pounds_{m_2,q}v_0\phi_2+h_2(x)|v_0|^{q-2}v_0\phi_2-F_{v_0}(x,u_0,v_0,w_0)\phi_2\right)d\mu=0,
\end{eqnarray*}
\end{small}which implies that $(u_0,v_0)$ solves (\ref{eq1}) for $w_0$.
\par
By $(u_n,v_n) \rightarrow (u_0,v_0)$ in $W$ and $C_1 \leq \|(u_n,v_n)\| \leq C_2$, as $n\to \infty$, we get
$$\|(u_0,v_0)\| \geq \|(u_n,v_n)\|-\|(u_n,v_n)-(u_0,v_0)\| \geq C_1-\|(u_n,v_n)-(u_0,v_0)\| \rightarrow C_1$$
and
$$\|(u_0,v_0)\| \leq \|(u_n,v_n)\|+\|(u_n,v_n)-(u_0,v_0)\| \leq C_2+\|(u_n,v_n)-(u_0,v_0)\| \rightarrow C_2.$$
Hence $C_1 \leq \|(u_0,v_0)\| \leq C_2$. \qed

\vskip2mm
\noindent
{\bf Remark 3.2.} There exist examples satisfying our results. For instance, let $ \gamma \in C(V)$ be such that $\min\limits_{x\in V}|\gamma(x)|>0$ and $w \in [-1,1]$, where $C(V)=\{u|u:V\to\R\}$.
Define
 $$
F(x,u,v,w)=(u^2(x)+v^2(x))^2(1+w^2)|\gamma(x)|.
$$
Then,
$$F_u(x,u,v,w)=4u(x)(u^2(x)+v^2(x))(1+w^2)|\gamma(x)|,$$
and
$$F_v(x,u,v,w)=4v(x)(u^2(x)+v^2(x))(1+w^2)|\gamma(x)|.$$
Set $p=3,\;q=2,\;\theta=r_1=r_2=4,\;c_1=c_2=16\|\gamma\|_\infty,\;a(|(u,v)|)=|(u,v)|^4=(u^2(x)+v^2(x))^2$ and $c(x)=|\gamma(x)|.$ Thus assumptions $(F_1)$, $(F_2)$ and $(H_1)$-$(H_3)$ are satisfied.\qed

\vskip2mm
{\section{Locally minimum type solutions}
\setcounter{equation}{0}
\vskip2mm
\par
In this section, we use the Ekeland's variational principle to find the local minimum solutions. We also present a dependence result on parameter  and a uniqueness result for the local minimum type solutions.
\vskip2mm
\noindent
{\bf Lemma 4.1.}  {\it Assume that $p=q$, $(F_1)$, $(F_2)$ and the following condition hold:\\
$(H_4)$  There exist $\delta>0$ and $L:V\rightarrow \R$ such that $L(x_0)>0$, $ \mu(x_0)L(x_0) >\frac{\|u^*\|_{W^{m_1,p}(V)}^p+\|v^*\|_{W^{m_2,p}(V)}^p}{p }$ for some $x_0\in V$ and $F(x_0,t,t,w)\geq L(x_0)t^p$ for all $0<t<\delta$, where $(u^*,v^*)$ is given in (\ref{EQ31}) below.\\
Then for all  $w\in \mathcal{J}$, $-\infty<\inf\{\varphi_w(u,v):(u,v)\in\bar B_{\rho} \}<0$, where $\rho$ is given in Lemma 3.2 in \cite{Zhang2022}, which is independent on $w$, and $\bar B_{\rho}=\{(u,v)\in W\big|\|(u,v)\|\leq \rho\}$.}
 \vskip2mm
\noindent
{\bf Proof.} Let
\begin{eqnarray}
\label{EQ31}
u^*(x)=v^*(x)=
\begin{cases}
1 \;\;\;\; \text{if}\;\; x=x_0,\\
0 \;\;\;\; \text{if}\;\; x \neq x_0.
\end{cases}
\end{eqnarray}
Then, for all  $w\in \mathcal{J}$ and all $t\in \R^+$, we have
\begin{eqnarray*}
\varphi_w(tu^*,tv^*)
&  =  & \frac{t^p}{p}\|u^*\|_{W^{m_1,p}(V)}^p+\frac{t^q}{q}\|v^*\|_{W^{m_2,q}(V)}^q - \int_VF(x,tu^*(x),tv^*(x),w)d\mu\nonumber\\
& \leq & \frac{t^p}{p}\|u^*\|_{W^{m_1,p}(V)}^p+\frac{t^q}{q}\|v^*\|_{W^{m_2,q}(V)}^q - \mu(x_0)F(x_0,t,t,w)\nonumber\\
& \leq & \frac{t^p}{p}\|u^*\|_{W^{m_1,p}(V)}^p+\frac{t^q}{q}\|v^*\|_{W^{m_2,q}(V)}^q - \mu(x_0) L(x_0)t^p.
\end{eqnarray*}
Note that $p>1$ and $p=q$. Then for all  $w\in \mathcal{J}$ and all $t\in \R^+$, by $(H_4)$, we have
\begin{eqnarray}
\label{EQ23}
\varphi_w(tu^*,tv^*) \leq  t^p \left(\frac{\|u^*\|_{W^{m_1,p}(V)}^p+\|v^*\|_{W^{m_2,p}(V)}^p}{p} - \mu(x_0) L(x_0)\right)<0.
\end{eqnarray}
Especially, we  choose $t_0$ satisfying
$$
 0<t_0<\min\left\{\delta,\frac{\rho}{\|(u^*,v^*)\|_W}\right\}
$$
such that $\varphi(t_0 u^*,t_0 v^*) < 0$. Obviously, $\|(t_0 u^*,t_0 v^*)\|_{W}<\rho$. Hence, $\inf\{\varphi(u,v):(u,v)\in\bar B_{\rho} \} \le \varphi(t_0 u^*,t_0 v^*)<0$. Moreover, from the proof of  Lemma 3.2 in \cite{Zhang2022}, it is easy to see that Lemma 3.2 in \cite{Zhang2022} is still true for $ \varphi_w$ and all $(u,v)\in\bar B_{\rho_w}$ and then by $(F_2)$, we have
\begin{eqnarray*}
       \varphi_w(u,v)
\geq    \left(\frac{1}{p}- \max\{K_1^pD, K_2^pD \} \right)\frac{\rho^p}{2^{p-1}},
\end{eqnarray*}
where $0<D<\min\{\frac{1}{pK^p_1},\frac{1}{pK^p_2}\}$. Hence,  $\varphi_w$ is bounded from below in $\bar B_{\rho_w}$ for all $w\in \mathcal{J}$. So $\inf\{\varphi_w(u,v):(u,v)\in\bar B_{\rho}\}>-\infty$. The proof is completed.
\qed

\vskip2mm
\noindent
{\bf Theorem 4.1} {\it Assume that $p=q$, $(F_1)$, $(F_2)$ and $(H_4)$ hold. Then, for any $w\in \mathcal{J}$, system (\ref{eq1}) has at least one solution $(u_w^*,v_w^*)\not=(0,0)$ of negative energy, which is locally minimum type, and
\begin{eqnarray}
\label{EQ120}
  C_3 \leq \|(u_w^*,v_w^*)\| \leq C_4,
\end{eqnarray}
where $C_3$ and $C_4$ are positive constants independent on $w$ and the values are given in $(\ref{EQ121})$ and $(\ref{EQ122})$ below,  respectively.}
\vskip2mm
\noindent
{\bf Proof.} Combining with Lemma 2.2 and Lemma 4.1 and being similar to the proof of Theorem 1.3 in \cite{Yang2024} with replacing $\varphi_\lambda, B_{\varrho_\lambda}$ with $\varphi_w, B_\rho$, respectively, there exists a sequence $\{(u_m,v_m)\}\subset \bar{B}_\rho $ such that
\begin{eqnarray}\label{aa1}
\varphi_w(u_m,v_m)\to \inf\limits_{\bar{B}_\rho}\varphi(u,v), \ \ \varphi_w'(u_m,v_m)\to 0, \ \ \mbox{as }m\to \infty.
\end{eqnarray}
Since $W$ is a finitely dimensional space, there exists a subsequence $\{(u_{m_i},v_{m_i})\}$ of $\{(u_m,v_m)\}\subset \bar{B}_\rho$ such that $(u_{m_i},v_{m_i})\to (u_w^*,v_w^*)$ in $W$ for some $(u_w^*,v_w^*) \in \bar{B}_\rho$. Then by (\ref{aa1}) and the continuity of $\varphi_w$ and  $\varphi'_w$,  $(u_w^*,v_w^*)$ is a critical point $\varphi_w$ and $\varphi(u_w^*,v_w^*)=\inf\limits_{\bar{B}_\rho}\varphi(u,v)$. By Lemma 4.1,  $\varphi(u_w^*,v_w^*)<0$, which, together with $(F_1)$,  implies that $(u_w^*,v_w^*)\not=(0,0)$.
\par
The proof of (\ref{EQ120}) is the same as (\ref{EQ4}) except for the value of $C_4$ which is mainly caused by  (\ref{aa2}). In details, $C_3$ is the same as (\ref{EQ17}) with $p=q$, that is,
\begin{small}
\begin{eqnarray}
\label{EQ121}
  C_3 =
\min\left\{\left({\frac{1}{2^p |V|\max{\{c_1,c_2\}}\max{\{b^{r_1},d^{r_2}\}}}}\right)^{\frac{1}{\max\{r_1,r_2\}-p}},
\left({\frac{1}{2^{p-1} |V|\max{\{c_1,c_2\}}\max{\{b^{r_1},d^{r_2}\}}}}\right)^{\frac{1}{\min\{r_1,r_2\}-p}}\right\}.
\end{eqnarray}
\end{small}We can obtain  $C_4$ just by replacing (\ref{aa2}) which is due to the mountain pass characterize by the following inequality which is due to the local minimizing characterize:
\begin{eqnarray*}
\label{EQ116}
\theta \varphi_w(u_w^*,v_w^*)
& \leq & \theta \varphi_w (t_0 u^*,t_0 v^*) \\
\notag
& \leq & \theta \left(\frac{t_0^p}{p}\|u^*\|_{W^{m_1,p}(V)}^p+\frac{t_0^p}{p}\|v^*\|_{W^{m_2,q}(V)}^p - \mu(x_0) L(x_0)t_0^p \right)\\
\notag
& \leq & \theta \left(\frac{t_0^p}{p}\|u^*\|_{W^{m_1,p}(V)}^p+\frac{t_0^p}{p}\|v^*\|_{W^{m_2,q}(V)}^p \right).
\end{eqnarray*}
Then combining the proof of Theorem 3.2, we can obtain that
\begin{small}
\begin{eqnarray}
\label{EQ122}
C_4 = \left({\frac{\theta 2^{p-1}t_0^p\left(\|u^*\|^p_{W^{m_1,p}(V)}+\| v^*\|^p_{W^{m_2,p}(V)}\right)}{p(\theta-p)}}\right)^{\frac{1}{p}}.
\end{eqnarray}
\end{small}
The proof is completed.\qed

\vskip2mm
\noindent
{\bf Theorem 4.2} {\it  Assume that $p=q$, $(F_1)$, $(F_2)$ and $(H_4)$  are satisfied. Let $\{w_n\} \subset \mathcal{J}$ be a convergent sequence of parameters with $\lim\limits_{n\rightarrow\infty}w_n=w_*$. For any nontrivial locally minimum type solutions sequence $\{(u_n,v_n)\}$ related to $\{w_n\}$ in (\ref{eq1}), there exists a subsequence $(u_{n_i},v_{n_i})\subset W$ such that $\lim\limits_{i\rightarrow\infty}(u_{n_i},v_{n_i})=(u_*,v_*)\in W$. Moreover, $(u_*,v_*)$ is a nontrivial solution to (\ref{eq1}) corresponding to $w_*$, and $C_3\leq\|(u_*,v_*)\|\leq C_4$.}
\vskip2mm
\noindent
{\bf Proof.}  The proof is the same as Theorem 3.3. We omit the details.\qed
\vskip2mm
\noindent
{\bf Remark 4.1.} In \cite{Galewski2011, Galewski2012, Galewski2013}, the authors mainly  considered the  parameter dependence of mountain pass type solutions, and they did not consider the parameter dependence of local minimum type solutions.  Theorem 4.1 and Theorem 4.2 is an attempt in this direction.
There is an interesting problem whether $(u_*,v_*)$ is still locally minimum type solutions (\ref{eq1}) corresponding to $w_*$ and whether $(u_0,v_0)$ is still the mountain pass type solution of (\ref{eq1}) corresponding to $w_0$. Currently, we have no idea to solve this problem.
\vskip2mm
\noindent
{\bf Remark 4.2.} We can get an example to support Theorem 4.1 and Theorem 4.2. For instance, let $ \gamma \in C(V)$ be such that $\min\limits_{x\in V}|\gamma(x)|>0$, $w \in [-1,1]$ and $0<e< \min\{\frac{1}{pK_1^p},\frac{1}{qK_2^q}\}$.
Define
$$
F(x,u,v,w)=e(u^2(x)+v^2(x))^2(1+w^2)|\gamma(x)|.
$$
Then,
$$F(x,u,u,w)=4eu^4(x)(1+w^2)|\gamma(x)|.$$
Set $p=q=4$ and $L(x_0)=4e \min\limits_{x\in V}|\gamma(x)|$ for some $x_0\in V$.
Thus assumptions $(F_1)$, $(F_2)$ and $(H_4)$ are satisfied.\qed
\vskip2mm
\par
Next, we present a uniqueness result for $(u^*_w,v^*_w)$. We need the following inequality.
\vskip2mm
\noindent
{\bf Lemma 4.2.}(\cite{Smejda2014}) {\it Assume that $p\geq2$, and there exists a positive constant $C_p$ such that
\begin{eqnarray*}
\label{EQ301}
(|x|^{p-2}x-|y|^{p-2}y)(x-y)\geq C_p|x-y|^p,\;\; \text{for all}\;\;x,y \in \R.
\end{eqnarray*}}
\par
 Let us introduce an additional condition:\\
{\it $(H_5)$ There exists two constants
$d_1,\;d_2 \in \left(0,\frac{C_p}{2^{p-1}|V|}\right)$
such that
$$|F_u(x,t_2,s_2,w)- F_u(x,t_1,s_1,w)|\leq d_1|t_2-t_1|^{p-1}$$
and
$$|F_v(x,t_2,s_2,w)- F_v(x,t_1,s_1,w)|\leq d_2|s_2-s_1|^{p-1}$$
 for all $x \in V$, $w \in \mathcal{J}$ and all $t_1,t_2, s_1,s_2 \in \R$ satisfying $|(t_1,t_2)|\le C_4\min\{\frac{1}{pK_1^p},\frac{1}{qK_2^q}\}$ and $|(s_1,s_2)|\le  C_4\min\{\frac{1}{pK_1^p},\frac{1}{qK_2^q}\}$.}
\vskip2mm
\noindent
{\bf Theorem 4.3.} {\it  Assume that $p=q$, $(F_1)$, $(F_2)$, $(H_4)$ and $(H_5)$  hold. Then for every fixed $w\in \mathcal{J}$, the solution $(u^*_w,v^*_w)$ of system (\ref{eq1}) in Theorem 4.1 is unique. }
\vskip2mm
\noindent
{\bf Proof.} Suppose that there exist two different functions $(u^*_{w1},v^*_{w1})$ and $(u^*_{w2},v^*_{w2})$ satisfying (\ref{eq1}). Then, by the proof of Lemma 3.4 with $p=q$ in \cite{Pang2023} and $(H_5)$, we can get
\begin{small}
\begin{eqnarray*}
\label{EQ302}
&&\frac{C_p}{2^{p-1}}\|(u^*_{w2}-u^*_{w1}, v^*_{w2}-v^*_{w1})\|^p  \\
&\leq & C_p\|u^*_{w2}-u^*_{w1}\|^p_{W^{m_1,p}(V)} + C_p\|v^*_{w2}-v^*_{w1}\|^p_{W^{m_2,p}(V)}    \\
&=& C_p\left(\int_V\left[\pounds_{m_1,p}(u^*_{w2}-u^*_{w1})(u^*_{w2}-u^*_{w1})+h_1(x)|u^*_{w2}-u^*_{w1}|^{p-2}(u^*_{w2}-u^*_{w1})(u^*_{w2}-u^*_{w1})\right]d\mu \right. \\
&& \left. +\int_V\left[\pounds_{m_2,p}(v^*_{w2}-v^*_{w1})(v^*_{w2}-v^*_{w1})+h_2(x)|v^*_{w2}-v^*_{w1}|^{p-2}(v^*_{w2}-v^*_{w1})(v^*_{w2}-v^*_{w1})\right]d\mu \right)\\
&\leq& \int_V\left[\pounds_{m_1,p}u^*_{w2}(u^*_{w2}-u^*_{w1})+h_1(x)|u^*_{w2}|^{p-2}u^*_{w2}(u^*_{w2}-u^*_{w1})\right]d\mu\\
&&- \int_V\left[\pounds_{m_1,p}u^*_{w1}(u^*_{w2}-u^*_{w1})+h_1(x)|u^*_{w1}|^{p-2}u^*_{w1}(u^*_{w2}-u^*_{w1})\right]d\mu\\
&&+\int_V\left[\pounds_{m_2,p}v^*_{w2}(v^*_{w2}-v^*_{w1})+h_2(x)|v^*_{w2}|^{p-2}v^*_{w2}(v^*_{w2}-v^*_{w1})\right]d\mu\\
&&- \int_V\left[\pounds_{m_2,p}u^*_{w1}(v^*_{w2}-v^*_{w1})+h_2(x)|v^*_{w1}|^{p-2}v^*_{w1}(v^*_{w2}-v^*_{w1})\right]d\mu\\
&=& \int_V [F_u(x,u^*_{w2},v^*_{w2},w)- F_u(x,u^*_{w1},v^*_{w1},w)](u^*_{w2}-u^*_{w1})d\mu + \int_V [F_v(x,u^*_{w2},v^*_{w2},w)- F_v(x,u^*_{w1},v^*_{w1},w)](v^*_{w2}-v^*_{w1})d\mu\\
&\leq& d_1|V|\|u^*_{w2}-u^*_{w1}\|^p + d_2|V|\|v^*_{w2}-v^*_{w1}\|^p\\
&\leq& \max \{d_1,d_2\}|V|\|(u^*_{w2}-u^*_{w1}, v^*_{w2}-v^*_{w1})\|^p.
\end{eqnarray*}
\end{small}
Thus,
$$\left(\frac{C_p}{2^{p-1}}- \max \{d_1, d_2\}|V|\right)\|(u^*_{w2}-u^*_{w1}, v^*_{w2}-v^*_{w1})\|^p\leq 0,$$
that is, $u^*_{w1}=u^*_{w2}$ and $v^*_{w1}=v^*_{w2}$, since $d_1, d_2<\frac{C_p}{2^{p-1}|V|}$. Therefore, $(u^*_{w1},v^*_{w1})=(u^*_{w2},v^*_{w2})$.\qed

\vskip2mm
\vskip2mm
\noindent
{\bf Remark 4.3.} We can present an example satisfying Theorem 4.3. Let $ \gamma \in C(V)$ be such that $\min\limits_{x\in V}|\gamma(x)|>0$, $w \in [-1,1]$ and $0<e^*< \min\{\frac{1}{pK_1^p},\frac{1}{qK_2^q}\}$.
Define
$$
F(x,u,v,w)=e^*(u^2(x)+v^2(x))(1+w^2)|\gamma(x)|.
$$
Then,
$$F_u(x,u,v,w)=2u(x)e^*(1+w^2)|\gamma(x)|.$$
and
$$F_v(x,u,v,w)=2v(x)e^*(1+w^2)|\gamma(x)|.$$
Set $p=q=2$, $d_1=d_2=4e^*\|\gamma\|_\infty$ and $L(x_0)=4e^* \min\limits_{x\in V}|\gamma(x)|$ for some $x\in V$.
Thus it is easy to verify that $(F_1)$, $(F_2)$, $(H_4)$ and $(H_5)$ are satisfied.\qed

\vskip2mm
{\section{An optimal control problem}}
\setcounter{equation}{0}
\par
In this section, we consider an optimal control problem. In details, we consider the minimization action functional
\begin{eqnarray*}
\psi(u,v,w)=\int_V g(x,u(x),v(x),w)d\mu,
\end{eqnarray*}
where $(u,v,w)$ satisfy (\ref{eq1}), and $g$ satisfies the following condition:\\
{\it $(G)$ $g: V\times \R^2 \times \mathcal{J} \rightarrow \R $ is continuous.}
\par
We fix a $w \in \mathcal{J}$ and  take $(u_w,v_w)$ as a solution to (\ref{eq1}) corresponding to $w$. We establish a set $A\subset W\times \mathcal{J} $,  where $A=\{(u_w,v_w,w)|(u_w,v_w,w) \mbox{ satisfies (\ref{eq1})} \}$.
\vskip2mm
\noindent
{\bf Remark 5.1.} Note that $\mathcal{J}$ is a bounded closed interval. Then for any $\{w_k\}\subset \mathcal{J}$, there exist a subsequence, still denoted by ${w_k}$, such that  $\lim\limits_{k\rightarrow\infty}w_k=\overline{w}\in \mathcal{J}$. Moreover, by Lemma 2.1 and Theorem 3.2 (or Theorem 4.1), there is some $l>0$ such that $\|(u_w,v_w)\|_\infty \leq l$ when $(u_w,v_w,w)\in A$.
\vskip2mm
\noindent
{\bf Theorem 5.1.} {\it Assume that conditions $(F_1)$, $(F_2)$, $(G)$ and $(H_1)$-$(H_3)$ hold. Let $\{w_k\} \subset \mathcal{J}$ be a convergent sequence with $\lim\limits_{k\rightarrow\infty}w_k=\overline{w}$. Then there exists a pair $(\overline{u}_{\overline{w}},\overline{v}_{\overline{w}},\overline{w})\in W\times \mathcal{J}$ such that $\psi(\overline{u}_{\overline{w}},\overline{v}_{\overline{w}},\overline{w})=\inf\limits_{(u_w,v_w,w)\in A}\psi(u_w,v_w,w).$}
\vskip2mm
\noindent
{\bf Proof.} Since $\|(u_w,v_w)\|_\infty \leq l$, $w \in \mathcal{J}$ which is a bounded closed interval and $g$ is continuous, the functional $\psi$ is bounded from below on $A$. Then there is a sequence $\{(u_{w_k},v_{w_k},w_k)\}\subset A$ such that
$$\lim\limits_{k\rightarrow\infty}\psi(u_{w_k},v_{w_k},w_k)=\inf\limits_{(u_w,v_w,w)\in A}\psi(u_w,v_w,w).$$
By Theorem 3.2 (or Theorem 4.1), $\{(u_{w_k},v_{w_k})\}$ is uniformly bounded in $W$. Then there exists a subsequence $\{(u_{w_{k_i}},v_{w_{k_i}},w_{k_i})\}$ such that $(u_{w_k},v_{w_k},w_{k_i})\rightarrow (\overline{u}_{\overline{w}},\overline{v}_{\overline{w}},\overline{w})$ in $W \times \mathcal{J}$. Then from Theorem 3.3 (or Theorem 4.2), we get $(\overline{u}_{\overline{w}},\overline{v}_{\overline{w}},\overline{w})$ is also a solution to (\ref{eq1}). Thus $(\overline{u}_{\overline{w}},\overline{v}_{\overline{w}},\overline{w}) \in A$. Therefore, according to the condition (G), it is easy to see that $\psi$ is continuous on $A$. Then, we can obtain
\begin{eqnarray*}
\inf\limits_{(u_w,v_w,w)\in A}\psi(u_w,v_w,w)
&  =  &  \lim\limits_{k\rightarrow\infty} \psi(u_{w_k},v_{w_k},w_k)\\
&  =  &  \lim\limits_{k\rightarrow\infty} \psi(u_{w_{k_i}},v_{w_{k_i}},w_{k_i})\\
&  =  &  \liminf\limits_{k\rightarrow\infty} \psi(u_{w_{k_i}},v_{w_{k_i}},w_{k_i})\\
& = & \psi(\overline{u}_{\overline{w}},\overline{v}_{\overline{w}},\overline{w}).
\end{eqnarray*}
Therefore, $(\overline{u}_{\overline{w}},\overline{v}_{\overline{w}},\overline{w})$ is a solution of the optimal control problem. The proof is completed. \qed
\vskip2mm
\noindent
{\bf Remark 5.2.} We present an example of $g$ which satisfies the assumptions of Theorem 5.1. Let
$$g(x,u,v,w)=z(x)(u^2+v^2)^2w^2,$$
where $z:V\to \R^+$.\qed

\vskip2mm
{\section{Nonexistence of solutions}
\setcounter{equation}{0}
\vskip2mm
\par
In this section, similar to \cite{Smejda2014}, we present a condition under which there is no solution to system (\ref{eq1}).
\vskip2mm
\noindent
{\bf Theorem 6.1.} {\it  Assume that $(F_1)$ holds and
\begin{eqnarray}
\label{EQ30}
 F_t(x,t,s,w)t+F_s(x,t,s,w)s <0 \;\; \text{for all} \;\; x\in V,\ (t,s)\in \R^2\backslash\{(0,0)\},w\in \mathcal{J}.
\end{eqnarray}
Then for every $w\in \mathcal{J}$, system (1.1) has no nontrivial solutions.}
\vskip2mm
\noindent
{\bf Proof.} Assume that (1.1) has a nontrivial solution $(\widetilde{u}_w,\widetilde{v}_w)$. Then $(\widetilde{u}_w,\widetilde{v}_w)$ is a nontrivial critical point of $\varphi_w$. According to (\ref{EQ6}), we get
\begin{small}
\begin{eqnarray*}
\label{EQ18}
0=\|\widetilde{u}_w\|^p_{W^{m_1,p}(V)}+\|\widetilde{v}_w\|^q_{W^{m_2,q}(V)}-\int_ VF_{\widetilde{u}_w}(x,\widetilde{u}_w(x),\widetilde{v}_w(x),w)\widetilde{u}_w(x) d\mu
-\int_ VF_{\widetilde{v}_w}(x,\widetilde{u}_w(x),\widetilde{v}_w(x),w)\widetilde{v}_w(x)d\mu.
\end{eqnarray*}
\end{small}Thus, it holds that
\begin{small}
\begin{eqnarray}
\label{EQ19}
&&\int_ VF_{\widetilde{u}_w}(x,\widetilde{u}_w(x),\widetilde{v}_w(x),w)\widetilde{u}_w(x) d\mu+\int_ VF_{\widetilde{v}_w}(x,\widetilde{u}_w(x),\widetilde{v}_w(x),w)\widetilde{v}_w(x)d\mu\nonumber\\
&=&\|\widetilde{u}_w\|^p_{W^{m_1,p}(V)}+\|\widetilde{v}_w\|^q_{W^{m_2,q}(V)} \geq 0
\end{eqnarray}
\end{small}
for all $w\in\mathcal{J}$. On the other hand, it follows from (\ref{EQ30}) that
\begin{eqnarray*}
\label{EQ20}
\int_ V (F_{\widetilde{u}_w}(x,\widetilde{u}_w(x),\widetilde{v}_w(x),w)\widetilde{u}_w(x)+ F_{\widetilde{v}_w}(x,\widetilde{u}_w(x),\widetilde{v}_w(x),w)\widetilde{v}_w(x))d\mu < 0
\end{eqnarray*}
for all $w\in\mathcal{J}$.  This contradicts with (\ref{EQ19}). So for every $w\in \mathcal{J}$, system (\ref{eq1}) has no nontrivial solutions.\qed
\vskip2mm
\noindent
{\bf Remark 6.1.} Let us take functions $F_u$: $V\times \R^2\times \mathcal{J}\rightarrow \R$ and $F_v$: $V\times \R^2\times \mathcal{J}\rightarrow \R$ given by
$$F_u(x,u,v,w)=-x^2 \arctan u,$$
and
$$F_v(x,u,v,w)=-x^2 \arctan v,$$
Then we see that assumptions of Theorem 6.1 are satisfied.

\vskip2mm
{\section{Results for the scalar equation (\ref{eq3})}}
\setcounter{equation}{0}
\vskip2mm

The results for system (\ref{eq1}) can also be applied to the following scalar equation (\ref{eq3}) on weighted finite graph:
\begin{eqnarray*}
\label{eq51}
  \pounds_{m,p}u+h(x)|u|^{p-2}u=f(x,u,w),\;\;\; x\in V,
\end{eqnarray*}
where $m\in \mathbb{N}$ and $p>1$, $h:V\to \R^+$ and $F(t)=\int_0^t f(x,t,w)dt,\; f:V\times \R\times \mathcal{J} \to \R$, $w\in \mathcal{J}$, where $\mathcal{J} \subset \R$ is a bounded closed interval.
\vskip2mm
\noindent
{\bf Theorem 7.1.} {\it Assume that $F$ satisfies the following conditions:\\
$(F'_1)$ \;  $F(x,0,w)=0$ and $F(x,t,w)$ is continuously differentiable in $ t\in \R$ for all $x\in V$, $w\in \mathcal{J}$;\\
$(F'_2)$ \; $\lim\limits_{|t|\rightarrow0}\frac{F(x,t,w)}{|t|^p}<\frac{1}{pK^P}$ for all $x\in V, w\in \mathcal{J}$, where
$
K=\frac{\left(\sum_{x\in V}\mu(x)\right)^{\frac{1}{p}}}{\mu_{\min}^\frac{1}{p}h_{\min}^\frac{1}{p}};
$\\
$(F'_3)$ \; $\lim\limits_{|t|\rightarrow\infty}\frac{F(x,t,w)}{|t|^p}=+\infty$ for all $x\in V, w\in \mathcal{J}$;\\
$(F'_4)$ \;  There exists a constant $\gamma_1>0$ such that
$$
\liminf\limits_{|t|\rightarrow\infty} \frac{f(x,t,w)t-pF(x,t,w)}{|t|^{\gamma_1}}
>0 \; \mbox{ for all} \;\; x\in V,w\in \mathcal{J}.
$$
Then equation (\ref{eq3}) has at least one nontrivial solution $u_w$.}
\vskip2mm
\noindent
{\bf Theorem 7.2.} {\it Assume that $F$ satisfies $(F'_1)$, $(F'_2)$ and the following conditions:\\
$(H'_1)$  $\theta F(x,t,w)\leq f(x,t,w)t$ for all  $x\in V, t\in \R \backslash\{0\},w\in \mathcal{J}$, where $\theta>p$;\\
$(H'_3)$  There exist two functions $a\in C(\R)$ and $c:V\rightarrow\R^+$ such that
$$F(x,t,w) \geq a(|t|)c(x) \; \mbox{ for all} \;\; x\in V, t\in \R,w\in \mathcal{J}.$$
Then  the equation (\ref{eq3}) has at least one nontrivial solution $u_w$.}
\vskip2mm
\noindent
{\bf Theorem 7.3.} {\it Assume that $F$ satisfies $(F'_1), (F'_2)$, $(H'_1)$, $(H'_3)$ and the following condition:\\
$(H'_2)$  $f(x,t,w)t\leq c_1 |t|^{r_1}$ for all  $x\in V, t\in \R,w\in \mathcal{J}$, where $c_1>0$ and $r_1>p$. \\
Then, for any $w\in \mathcal{J}$, the equation (\ref{eq3}) has at least one nontrivial solution $u_w$ satisfying
\begin{eqnarray*}
\label{EQq}
  C'_1 \leq \|u_w\| \leq C'_2,
\end{eqnarray*}
where $C'_1$ and $C'_2$ are positive constants independent on $w$ and the values are given in (\ref{EQ112}) and (\ref{EQ113}) below, respectively.}
\vskip2mm
\noindent
{\bf Proof.} Choose arbitrary but fixed $w\in \mathcal{J}$. According to Theorem 7.1, we know that the equation (\ref{eq3}) has  one nontrivial solution $u_w$. Then from the proof of \cite{Zhang2022}, there exists a $0\neq u_{0}\in \R$  which is independent on $w$ since $(F_3)'$ holds for all $w\in \mathcal{J}$ (see \cite{Zhang2022}), such that
\begin{eqnarray*}
\label{EQ111}
  \varphi_w(u_w)=\inf_{\gamma\in\Gamma} \max_{t\in[0,1]} \varphi_w(\gamma(t)):=C'^*>0,
\end{eqnarray*}
where
$$
\Gamma:=\{\gamma\in C([0,1],\R):\gamma(0)=0,\gamma(1)=u_{0}\}.
$$
Similar to the proof of Theorem 3.2,  we can get
\begin{eqnarray}
\label{EQ112}
  \|u_w\| \geq C'_1=\left(\frac{1}{2^p|V|c_1 b^{r_1}}\right)^\frac{1}{r_1-p},
\end{eqnarray}
and
\begin{eqnarray}
\label{EQ113}
  \|u_w\| \leq C'_2=\left(\frac{\theta 2^{p-1} \|u_0\|^p}{\theta-p}\right)^{\frac{1}{p}}.
\end{eqnarray}
The proof is completed.

\vskip2mm
\noindent
{\bf Theorem 7.4.} {\it Assume that $(F'_1)$, $(F'_2)$ and $(H'_1)$-$(H'_3)$ are satisfied. Let $\{w_n\} \subset \mathcal{J}$ be a convergent sequence of parameters with $\lim\limits_{n\rightarrow\infty}w_n=w_0$. For any mountain pass type nontrivial solutions sequence $\{u_n\}$ related to ${w_n}$ in (\ref{eq3}), which are given in Theorem 7.3, there exist a subsequence $
\{u_{n_i}\} \subset \R$ such that $\lim\limits_{i\rightarrow\infty}u_{n_i}=u_0\in \R$. Moreover, $u_0$ is a nontrivial solution to (\ref{eq3}) corresponding to $w_0$, and $C'_1\leq\|u_0\|\leq C'_2$.}

\vskip2mm
\noindent
{\bf Theorem 7.5.}  {\it Assume that $F$ satisfies $(F'_1)$, $(F'_2)$ and the following condition:\\
$(H'_4)$  There exist $\delta'>0$ and  $L:V\rightarrow \R$ such that $L(x_0)>0$, $ \mu(x_0)L(x_0) >\frac{\|u^*\|^p}{p}$ for some $x_0\in V$ and $F(x,t,w)\geq L(x_0)t^p$ for all $0<t<\delta'$, where $u^*$ is given by (\ref{EQ31}).\\
Then, any $w\in \mathcal{J}$, system (\ref{eq3}) has at least one solution $u_w^*\not=0$ of negative energy, which is locally minimum type, and
\begin{eqnarray}
\label{EQ200}
  C'_3 \leq \|u_w^*\| \leq C'_4,
\end{eqnarray}
where $C'_3=C'_1$ and $C'_4=\left(\frac{\theta 2^{p-1}t_0^p\|u^*\|^p}{\theta-p}\right)^{\frac{1}{p}}$}.

\vskip2mm
\noindent
{\bf Theorem 7.6} {\it  Assume that $(F'_1)$, $(F'_2)$ and $(H'_4)$  are satisfied. Let $\{w_n\} \subset \mathcal{J}$ be a convergent sequence of parameters with $\lim\limits_{n\rightarrow\infty}w_n=w_*$. For any locally minimum tuype nontrivial solutions sequence $u_n$ related to ${w_n}$ in (\ref{eq3}), which are given in Theom 7.4, there exists a subsequence $u_{n_i}\subset \R$ such that $\lim\limits_{i\rightarrow\infty}u_{n_i}=u_*\in \R$. Moreover, $u_*$ is a nontrivial solution to (\ref{eq3}) corresponding to $w_*$, and $C'_3\leq\|u_*\|\leq C'_4$.}

\vskip2mm
\noindent
{\bf Theorem 7.7.} {\it Assume that $(F'_1)$, $(F'_2)$, $(H'_1)$-$(H'_3)$ and the following condition hold:\\
$(G')$ $g: V\times \R \times \mathcal{J} \rightarrow \R $ is continuous and for any fixed $(x,u)\in V\times \R $, the functional $w\rightarrow g(x,u,w)$ is convex.\\
   For any fixed  $w \in \mathcal{J}$,  take $u_w$ as a solution to (\ref{eq3}) corresponding to $w$. Let $A'=\{(u_w,w)|(u_w,w) \mbox{ is a solution of (\ref{eq3})}\}$ and
$\{w_k\} \subset \mathcal{J}$ be a convergent sequence with $\lim\limits_{k\rightarrow\infty}w_k=\overline{w}$. Then there exists a pair $(\overline{u}_{\overline{w}},\overline{w})\in \R \times \mathcal{J}$ such that $\psi(\overline{u}_{\overline{w}},\overline{w})=\inf\limits_{(u_w,w)\in A'}\psi(u_w,w).$}

\vskip2mm
\noindent
{\bf Theorem 7.8.} {\it Let $(F'_1)$ hold. Assume that
\begin{eqnarray}
\label{EQ201}
 F_t(x,t,s,w)t <0 \;\; \text{for all} \;\; x\in V,\ t\in \R\backslash\{0\},w\in \mathcal{J}.
\end{eqnarray}
Then for every $w\in \mathcal{J}$, the problem (\ref{eq3}) has no nontrivial solutions.}

 \vskip2mm
 \noindent
 {\bf Acknowledgments}\\
This work is supported by Yunnan Fundamental Research Projects of China (grant No: 202301AT070465) and  Xingdian Talent
Support Program for Young Talents of Yunnan Province in China.

 \vskip2mm
 \noindent
 {\bf Competing interests}\\
The authors declare that they have no competing interests.

\vskip2mm
\renewcommand\refname{References}
{}

\begin{thebibliography}{}

\bibitem{Adimurthi2010}A. Adimurthi, Y. Yang. An interpolation of Hardy inequality and Trudinger-Moser inequality in $\mathbb{R}^N$ and its applications. Int. Math. Res. Notices, 2010(13), 2394-2426.

\bibitem{Ambrosetti1973}A. Ambrosetti, P. Rabinowitz. Dual variational methods in critical point theory and applications. J. Funct. Anal, 1973, 14, 349-381.

\bibitem{Molica Bisci2017}G. Bisci, D. Repov\v{s}, R. Servadei. Nonlinear problems on the Sierpi\'{n}ski gasket. J. Math. Anal. Appl, 2017, 452, 883-895.

\bibitem{Bonanno2012}G. Bonanno, G. Bisci, V. R\v{a}dulescu. Variational analysis for a nonlinear elliptic problem on the Sierpi\'{n}ski gasket. Esaim Contr. Optim. Ca, 2012, 18(4), 941-953.

\bibitem{Breckner2010}B. Breckner, D. Repov\v{s}, C. Varga. On the existence of three solutions for the Dirichlet problem on the Sierpi\'{n}ski gasket. Nonlinear Anal-Theor, 2010, 73(9), 2980-2990.

\bibitem{Breckner2011}B. Breckner, V. R\v{a}dulescu, C. Varga. Infinitely many solutions for the Dirichlet problem on the Sierpi\'{n}ski gasket. Anal. Appl, 2011, 9(3), 235-248.

\bibitem{Chang1986}K. Chang. Critical point theory and its applications (in Chinese), Xiandai Shuxue Congshu. Shanghai Kexue Jishu Chubanshe, Shanghai, 1986.

\bibitem{Falconer1999}K. Falconer, J. Hu. Nonlinear elliptical equations on the Sierpi\'{n}ski gasket. J. Math. Anal. Appl, 1999, 240, 552-573.

\bibitem{Galewski2019}M. Galewski. On the mountain pass solution to boundary value problems on the Sierpi\'{n}ski gasket. Results Math, 2019, 74, 167.

\bibitem{Galewski2013}M. Galewski, J. Smejda. On the dependence on parameters for mountain pass solutions of second order discrete BVP's. Appl. Math. Comput, 2013, 219(11), 5963-5971.

\bibitem{Galewski2011}M. Galewski. Dependence on parameters for a discrete Emden-Fowler equation. Appl. Math. Comput, 2011, 218(4), 1247-1253.

\bibitem{Galewski2012}M. Galewski. A note on the dependence on parameters for a nonlinear system via monotonicity theory. J. Differ. Equ. Appl, 2012, 18(7), 1253-1256.

\bibitem{Ge20181}H. Ge. A $p$-th Yamabe equation on graph. P. Am. Math. Soc, 2018, 146(5), 2219-2224.

\bibitem{Ge20184}H. Ge, W. Jiang. Yamabe equations on infinite graphs. J. Math. Anal. Appl, 2018, 460(2), 885-890.

\bibitem{Grigor yan2016}A. Grigor¡¯yan, Y. Lin, Y. Yang. Yamabe type equations on graphs. J. Differ. Equations, 2016, 261(9), 4924-4943.

\bibitem{Han2021}X. Han, M. Shao. $p$-Laplacian equations on locally finite graphs. Acta Math. Sin, 2021, 37(11), 1-34.

\bibitem{Mawhin1989}J. Mawhin, M. Willem. Critical point theorem and Hamiltonian seytems. Appl. Math. Sci. Springer-Verlag, New York, 1989, 74.

\bibitem{Pang2023}Y. Pang, X. Zhang. Existence of three solutions for a poly-Laplacian system on graphs. arXiv preprint arXiv: 2309.09849, 2023.

\bibitem{Pang2024}Y. Pang, J. Xie, X. Zhang. Infinitely many solutions for three quasilinear Laplacian systems on weighted graphs. Bound. Value Probl, 2024, 45, 1-23.

\bibitem{Rabinowitz1986}P. Rabinowitz. Minimax methods in critical point theory with applications to differential equations. J. Am. Math. Soc, 1986.

\bibitem{Smejda2014} J. Smejda, R. Wieteska. On the dependence on parameters for second order discrete boundary value problems with the $p(k)$-Laplacian. J. Opuscu. Math, 2014, 34(4), 851-870.

\bibitem{Yang2024}P. Yang, X. Zhang. Existence and multiplicity of nontrivial solutions for a $(p,q)$-Laplacian system on locally finite graphs. Taiwanese. J. Math. 2024, 28 (3), 551-588

\bibitem{Yang2023}P. Yang, X. Zhang. Existence of solutions for a poly-Laplacian system involving concave-convex nonlinearity on locally finite graphs. Electron. Res. Arch, 2023, 31(12), 7473-7495.

\bibitem{Yu2024}X. Yu, X. Zhang, J. Xie, X. Zhang. Existence of nontrivial solutions for a class of poly-Laplacian system with mixed nonlinearity on graphs. Math. Methods Appl. Sci, 2024, 47(4), 1750-1763.

\bibitem{Zhang2018}X. Zhang, A. Lin. Positive solutions of $p$-th Yamabe type equations on graphs. Front. Math. China, 2018, 13(6), 1501-1514.

\bibitem{Zhang2019}X. Zhang, A. Lin. Positive solutions of $p$-th Yamabe type equations on infinite graphs. P. Am. Math. Soc, 2019, 147(4), 1421-1427.

\bibitem{Zhang2012}X. Zhang, X. Tang. Non-constant periodic solutions for second order Hamiltonian system with a p-Laplacian. J. Math. Slovaca, 2012, 62(2), 231-246.

\bibitem{Zhang2022}X. Zhang, X. Zhang, J. Xie, X. Yu. Existence and multiplicity of nontrivial solutions for poly-Laplacian systems on finite graphs. Bound. Value Probl, 2022, 2022(1), 1-13.













\end{thebibliography}
\end{document}